\documentclass{article}
\usepackage{CJK,CJKnumb,CJKulem,times,dsfont,ifthen,mathrsfs,latexsym,amsfonts,color}
\usepackage{amsmath,amsthm,makeidx,fontenc,amssymb,bm,graphicx,psfrag,listings,curves,extarrows}
\usepackage[colorinlistoftodos]{todonotes}
\makeindex
\newtheorem{theorem}{Theorem}[section]
\newtheorem{definition}[theorem]{Definition}
\newtheorem{lemma}[theorem]{Lemma}
\newtheorem{corollary}[theorem]{Corollary}
\newtheorem{proposition}[theorem]{Proposition}
\newtheorem{remark}[theorem]{Remark}

\newcommand{\s}{\section}

\newcommand{\la}{\lambda}

\newcommand{\pa}{\partial}

\newcommand{\R}{\mathbb R}
\newcommand{\al}{\alpha}

\newcommand{\C}{\mathbb C}

\newcommand{\lab}{\label}
\newcommand{\bt}{\begin{theorem}}
\newcommand{\et}{\end{theorem}}
\newcommand{\bl}{\begin{lemma}}
\newcommand{\el}{\end{lemma}}
\newcommand{\bd}{\begin{definition}}
\newcommand{\ed}{\end{definition}}
\newcommand{\bc}{\begin{corollary}}
\newcommand{\ec}{\end{corollary}}
\newcommand{\bp}{\begin{proof}}
\newcommand{\ep}{\end{proof}}
\newcommand{\bx}{\begin{example}}
\newcommand{\ex}{\end{example}}
\newcommand{\bi}{\begin{exercise}}
\newcommand{\ei}{\end{exercise}}
\newcommand{\bo}{\begin{proposition}}
\newcommand{\eo}{\end{proposition}}
\newcommand{\br}{\begin{remark}}
\newcommand{\er}{\end{remark}}
\newcommand{\beq}{\begin{equation}}
\newcommand{\eeq}{\end{equation}}
\newcommand{\ba}{\begin{align}}
\newcommand{\ea}{\end{align}}
\newcommand{\bn}{\begin{enumerate}}
\newcommand{\en}{\end{enumerate}}
\newcommand{\bg}{\begin{align*}}
\newcommand{\bcs}{\begin{cases}}
\newcommand{\ecs}{\end{cases}}

\newcommand{\ga}{\gamma}

\newcommand{\bean}{\begin{eqnarray*}}
\newcommand{\eean}{\end{eqnarray*}}

\newenvironment{altproof}[1]
{\noindent%\addvspace{0.3cm}
{\sl Proof of {#1}}.}
{\nopagebreak\mbox{}\hfill $\Box$\par\addvspace{0.5cm}}

\newcommand{\ind}{\text{\rm ind}}
\def\A{\mathbb{A}}
\def\C{\mathbb{C}}
\def\N{\mathbb{N}}

\def\P{\mathbb{P}}
\def\R{\mathbb{R}}

\def\bd{\mathrm{bd}\,}

\newcommand{\cb}{{\mathcal B}}
\newcommand{\cc}{{\mathcal C}}
\newcommand{\cd}{{\mathcal D}}

\newcommand{\co}{{\mathcal O}}

\newcommand{\cs}{{\mathcal S}}
\newcommand{\ct}{{\mathcal T}}

\newcommand{\cz}{{\mathcal Z}}

\newcommand{\eps}{\varepsilon}
\newcommand{\be}{\beta}
\newcommand{\de}{\delta}
\newcommand{\ka}{\kappa}
\newcommand{\om}{\omega}

\newcommand{\De}{\Delta}

\newcommand{\Om}{\Omega}

\newcommand{\weakto}{\rightharpoonup}

\numberwithin{equation}{section}
\begin{document}
\begin{CJK*}{GBK}{song}
\title{\bf  Normalized solutions for a coupled Schr\"odinger system\thanks{Supported by NSFC(11801581), Guangdong NSFC(2018A030310082),
E-mails: thomas.bartsch@math.uni-giessen.de(Bartsch)\quad\quad\quad zhongxuexiu1989@163.com(Zhong)\quad\quad\quad  zou-wm@mail.tsinghua.edu.cn(Zou)}}

\date{}
\author{
{\bf Thomas Bartsch$^1$,\; \;Xuexiu Zhong$^2$\;\&\;Wenming Zou$^3$}\\
\footnotesize \it 1. Mathematisches Institut, Justus-Liebig-Universit\"at Giessen, Arndtstrasse 2, 35392 Giessen, Germany.\\
\footnotesize \it 2. South China Research Center for Applied Mathematics and Interdisciplinary Studies, \\
\footnotesize \it South China Normal University, Guangzhou 510631, China.\\
\footnotesize \it 3. Department of Mathematical Sciences, Tsinghua University, Beijing 100084, China.}

\maketitle
\begin{center}
\begin{minipage}{120mm}
\begin{center}{\bf Abstract}\end{center}
In the present paper, we prove the existence of solutions $(\la_1,\la_2,u,v)\in\R^2\times H^1(\R^3,\R^2)$ to systems of coupled Schr\"odinger equations
$$
\bcs
-\De u+\la_1u=\mu_1 u^3+\be uv^2\quad &\hbox{in}\;\R^3\\
-\De v+\la_2v=\mu_2 v^3+\be u^2v\quad&\hbox{in}\;\R^3\\
u,v>0&\hbox{in}\;\R^3
\ecs
$$
satisfying the normalization constraint
$
\displaystyle\int_{\R^3}u^2=a^2\quad\hbox{and}\;\int_{\R^3}v^2=b^2,
$
which appear in binary mixtures of Bose-Einstein condensates or in nonlinear optics.
The parameters $\mu_1,\mu_2,\be>0$ are prescribed as are the masses $a,b>0$. The system has been considered mostly in the fixed frequency case. And  when the masses are prescribed, the standard approach to this problem is variational with $\la_1,\la_2$ appearing as Lagrange multipliers. Here we present a new approach based on bifurcation theory and the continuation method. We  obtain the existence of normalized solutions for any given $a,b>0$ for $\be$ in a large range. We also give a result about the nonexistence of positive solutions. From which one can see that our existence theorem is almost the best. Especially, if $\mu_1=\mu_2$ we prove that normalized solutions exist for all $\be>0$ and all $a,b>0$.
 \vskip0.13in

{\it Key  words:}   Schr\"odinger system; self-focusing; attractive interaction; solitary wave; normalized solution; global bifurcation.


\vskip0.13in

\end{minipage}
\end{center}
\vskip0.26in
\newpage
%%%
\s{Introduction}
%%%
The time-dependent system of coupled nonlinear Schr\"odinger equations
\beq\lab{eq:NLS}
\begin{cases}
-i\frac{\partial}{\partial t}\Phi_1=\De \Phi_1+\mu_1 |\Phi_1|^2\Phi_1+\be |\Phi_2|^2\Phi_1,\\
-i\frac{\partial}{\partial t}\Phi_2=\De \Phi_2+\mu_2|\Phi_2|^2\Phi_1+\be |\Phi_1|^2\Phi_2,\\
\Phi_j=\Phi_j(x,t)\in \C, j=1,2,\;N\leq 3,
\end{cases}\;\;(x,t)\in \R^N\times \R,
\eeq
is used as model for various physical phenomena, for instance binary mixtures of Bose-Einstein condensates, or the propagation of mutually incoherent wave packets in nonlinear optics; see e.g.\ \cite{AkhmedievAnkiewicz.1999, Esry1998, Frantzeskakis2010, Timmermans1998}. In the models, $i$ is the imaginary unit, $\Phi_j$ is the wave function of the $j$-th component, and the real numbers $\mu_j$ and $\be$ represent the intra-spaces and inter-species scattering length, describing respectively the interaction between particles of the same component or of different components. In particular, the positive sign of $\mu_j$ (and of $\be$) stays for attractive interaction, while the negative sign stays for repulsive interaction. In present paper, we consider the case of positive parameters $\mu_1,\mu_2,\be>0$, i.e.\ the self-focusing and attractive case. An important, and of course well known, feature of \eqref{eq:NLS} is conservation of masses: the $L^2$-norms $|\Phi_1(\cdot,t)|_2, |\Phi_2(\cdot,t)|_2$ of solutions are independent of $t\in\R$. These norms have a clear physical meaning. In the aforementioned contexts, they represent the number of particles of each component in Bose-Einstein condensates, or the power supply in the nonlinear optics framework.

The ansatz $\Phi_1(x,t)=e^{i\la_1 t}u(x)$ and $\Phi_2(x,t)=e^{i\la_2t}v(x)$ for solitary wave solutions leads to the  elliptic system:
\beq\lab{eq:system}
\begin{cases}
  -\De u+\la_1 u=\mu_1 u^3+\be uv^2,\\
  -\De v+\la_2 v=\mu_2 v^3+\be vu^2,
\end{cases}\;\hbox{in}\;\R^N.
\eeq
This system has been investigated by many authors since about 2005, mainly in the fixed frequency case where $\la_1,\la_2>0$ are prescribed; see e.g.\ \cite{Ambrosetti2007, BartschWang2006, BartschWangWei2007, Chen2013a, LinWei.2005, MaiaMontefuscoPellacci.2006, Mandel2016, Sirakov2007, Soave2015, Soave2016, Terracini2009,Wei2008} and the references therein.

Much less is known when the $L^2$-norms $|u|_2,|v|_2$ are prescribed, in spite of the physical relevance of normalized solutions. A natural approach to finding solutions of \eqref{eq:system} satisfying the normalization constraints
\beq\lab{eq:norm}
  \int_{\R^N}u^2=a^2\quad\hbox{and}\;\int_{\R^N}v^2=b^2
\eeq
consists in finding critical points $(u,v)\in H^1(\R^N,\R^2)$ of the energy
\[
  J(u,v) = \frac12\int_{\R^N}\left(|\nabla u|^2+|\nabla v|^2\right)
               - \frac14\int_{\R^N}\left(\mu_1u^4+\mu_2v^4+2\be u^2v^2\right)
\]
under the constraints \eqref{eq:norm}. Then the parameters $\la_1,\la_2$ appear as Lagrange multipliers. All papers on normalized solutions of \eqref{eq:system} are based on this approach; see \cite{Bartsch2018, Bartsch2016, Bartsch2017, Bartsch2019, GouJeanjean2018} and the references therein. Only the papers \cite{Bartsch2016, GouJeanjean2018} deal with \eqref{eq:system}-\eqref{eq:norm} with $\be>0$. The existence of normalized solutions for systems of nonlinear Schr\"odinger equations with trapping potential has been proved in \cite{Noris2015}, and on bounded domains in \cite{Noris2019}, also by variational methods. In \cite{Noris2015,Noris2019} the masses $a^2,b^2$ have to be small.

In the present paper we propose a different approach based on bifurcation theory applied to \eqref{eq:system} with $\la_2=1$, taking $\la_1$ as parameter. There are two families of semitrivial solutions of \eqref{eq:system} where either $u=0$ or $v=0$. The bifurcation of global continua of positive solutions of \eqref{eq:system} from these semitrivial solutions has been proved in \cite{BartschWangWei2007}. We shall investigate the global behavior of these continua, and the $L^2$-norms of the solutions along them, in order to obtain the existence of solutions of \eqref{eq:system}-\eqref{eq:norm}. A major tool will be the fixed point index in cones.

In this paper we deal with the case $N=3$ when the growth of the nonlinearity is mass-supercritical. In dimension $N=1$ the growth of the nonlinearity is mass-subcritical so that $J$ is bounded from below on the constraint and normalized solutions can be obtained by minimization. In dimension $N=2$ the growth of the nonlinearity in \eqref{eq:system} is mass-critical making the existence of normalized solutions a very subtle issue, heavily depending on the prescribed masses $a^2,b^2$, as can already be seen in the scalar case.

The paper is organized as follows. In the next section we state and discuss our results, in particular we compare them with existing results on normalized solutions. We also state and discuss some new non-existence and uniqueness theorems for \eqref{eq:system} that will enter in the proofs of our results on normalized solutions. Then in Section~\ref{sec:prelim} we collect and prove a few basic facts about \eqref{eq:system}. Section~\ref{sec:global} contains the main idea of our approach. There we reduce the proofs of our results on normalized solutions to the problem of controlling the $L^2$-norms along continua of solutions of \eqref{eq:system}, and we describe the bifurcating continua. An important part of our proof is to understand the behavior of the $L^2$-norms as $\la\to0$ or $\la\to\infty$. We investigate this in Section~\ref{sec:asympt} where we also prove the non-existence and uniqueness theorems for \eqref{eq:system}. The main results about normalized solutions will be proved in section~\ref{sec:proofs}.

%%%
\s{Statement of results}
%%%
We are concerned with the existence of real numbers $\la_1,\la_2\in \R$ and of radial functions $u, v\in H_{rad}^1(\R^3)$ that solve
\beq\lab{eq:main}
\begin{cases}
-\De u+\la_1u=\mu_1 u^3+\be uv^2,\quad &\hbox{in}\;\R^3,\\
-\De v+\la_2v=\mu_2 v^3+\be u^2v,\quad&\hbox{in}\;\R^3,\\
u,v>0,&\hbox{in}\;\R^3,\\
|u|_2=a \quad\text{and}\quad |v|_2=b,
%\int_{\R^3}u^2=a^2\quad\hbox{and}\;\int_{\R^3}v^2=b^2,
\end{cases}
\eeq
where $\mu_1, \mu_2, \be, a, b>0$ are prescribed positive real numbers. In order to state our results we define
\beq\lab{eq:def-tau0}
  \tau_0:=\inf_{\phi\in \cd^{1,2}_0(\R^N)\setminus\{0\}}\frac{\int_{\R^3}|\nabla \phi|^2dx}{\int_{\R^3}U^2\phi^2dx},
  \eeq
where $U$ is the unique positive radial solution to
\beq\lab{eq:def-U}
  -\De u+u=u^3\;\hbox{in}\;\R^N;\quad u(x)\to0\ \text{ as $|x|\to\infty$;}
\eeq
cf.\ \cite{Kwong.1989}. We shall see that $\tau_0\in (0,1)$.

\bt\lab{Main-th}
Let $\mu_1,\mu_2>0$. Then we have the following conclusions.
\begin{itemize}
\item[a)] If $\be\in (0,\tau_0\min\{\mu_1,\mu_2\}]\cup (\tau_0\max\{\mu_1,\mu_2\},+\infty)$ then for any $a,b>0$, the  problem \eqref{eq:main} has a solution $(\la_1,\la_2,u,v)$ with $\la_1>0,\la_2>0$ and $u,v\in H_{rad}^1(\R^3)$.
\item[b)] If $\be\in (\tau_0\min\{\mu_1,\mu_2\},\tau_0\max\{\mu_1,\mu_2\}]$ then there exists $\de>0$ such that for any $a,b>0$ satisfying
$$
\bcs
   \frac{a}{b}\le \delta\quad &\hbox{if}\;\mu_2<\mu_1;\\
   \frac{a}{b}\ge \frac{1}{\delta}\quad &\hbox{if}\;\mu_2> \mu_1,
\ecs
$$
the problem \eqref{eq:main} has a solution $(\la_1,\la_2,u,v)$ with $\la_1>0,\la_2>0$ and $u,v\in H_{rad}^1(\R^3)$. If in addition $\be \in (\tau_0\min\{\mu_1,\mu_2\},\min\{\mu_1,\mu_2\})$ then
\[
  \de \ge \sqrt{\frac{\be-\min\{\mu_1,\mu_2\}}{\be-\max\{\mu_1,\mu_2\}}}.
\]
\end{itemize}
\et

%In the symmetric case $\mu_1=\mu_2$ we have a stronger result.

%\bt\lab{thm:mu1=mu2}
%If $\mu_1=\mu_2>0$, then for any prescribed $a,b>0$ and $\be>0$, the problem \eqref{eq:main} has a solution $(\la_1,\la_2,u,v)$ with $\la_1>0,\la_2>0$ and $u,v\in H_{rad}^1(\R^3)$.
%\et

Of course it is natural to ask whether \eqref{eq:main} has a solution without any conditions on $\mu_1,\mu_2,\beta,a,b$. This is not true however, as the next result shows.

\bo\lab{prop:nonexist-1}
If $\mu_2\le \be\le \tau_0\mu_1$, then there exists $q >0 $ such that \eqref{eq:main} has no solution for $\frac{a}{b} > q$. If $\mu_1\le \be\le \tau_0\mu_2$, then there exists $\tilde{q}>0$ such that \eqref{eq:main} has no solution for $\frac{a}{b} < \tilde{q}$.
\eo

Theorem~\ref{Main-th} and Proposition~\ref{prop:nonexist-1} will be proved in Section~\ref{sec:proofs}.

\br\lab{rem:compare}
As mentioned in the introduction, only the papers \cite{Bartsch2016, GouJeanjean2018} deal with \eqref{eq:system}-\eqref{eq:norm} in the case $\be>0$. Theorem~\ref{Main-th} significantly improves and complements the results of \cite{Bartsch2016}. There the authors obtain a solution $(\la_1,\la_2,u,v)$ of \eqref{eq:main} as in Theorem~\ref{Main-th} for $0<\be<\be_1$ and for $\be>\be_2$ where $\be_1,\be_2>0$ are defined implicitely by
%\beq\lab{20180720-e1}
\[
  \max\left\{\frac{1}{a^2\mu_1^2},\frac{1}{b^2\mu_2^2}\right\}=\frac{1}{a^2(\mu_1+\be_1)^2}+\frac{1}{b^2(\mu_2+\be_1)^2}.
\]
%\eeq
and
%\beq\lab{20180720-e2}
\[
  \frac{(a^2+b^2)^3}{(\mu_1a^4+\mu_2b^4+2\be_2a^2b^2)^2}=\min\left\{\frac{1}{a^2\mu_1^2},\frac{1}{b^2\mu_2^2}\right\}.
\]
%\eeq
Clearly the bounds $\be_1,\be_2$ depend on the masses $a,b>0$ and
\[
  \be_1\to0,\ \be_2\to\infty \quad\text{as $\frac{a}{b}\to0$ or $\frac{a}{b}\to\infty$}.
\]
In particular there is no value of $\be$ so that the results from \cite{Bartsch2016} yield a solution for all masses.\\
In \cite{GouJeanjean2018} the authors consider more general (but still homogeneous) nonlinearities and interaction terms. Specialized to \eqref{eq:system}-\eqref{eq:norm} their results recover those of \cite{Bartsch2016}. Our new approach via bifurcation theory and continuaton can also be applied to the systems considered in \cite{GouJeanjean2018} and to improve the results in that paper.
\er

We now add a few results on \eqref{eq:system} which enter in the proofs of Theorem~\ref{Main-th} and which have some interest in itself. Below we assume $\la_1,\la_2>0$. This is no restriction because we shall prove that positive solutions of \eqref{eq:system} with $\mu_1,\mu_2,\be>0$ can only exist if $\la_1,\la_2>0$; see Lemma~\ref{lem:sol-properties}.

\bt\lab{thm:nonexist-2}
\begin{itemize}
\item[\rm a)] For $\be\ge\mu_1$ there exists $\eta_1(\be)>0$ such that \eqref{eq:system} has no positive solution if $\frac{\la_1}{\la_2}>\eta_1(\be)$.
\item[\rm b)] For $\be\ge\mu_2$ there exists $\eta_2(\be)>0$ such that \eqref{eq:system} has no positive solution if $\frac{\la_1}{\la_2}<\eta_2(\be)$.
\end{itemize}
\et

The next theorem makes some progress towards uniqueness of positive solutions of \eqref{eq:system}.

\bt\lab{thm:sirakov}
\begin{itemize}
\item[\rm a)] Problem~\eqref{eq:system} with $N=3$ has at most one positive solution for $\frac{\la_1}{\la_2}>0$ small or for $\frac{\la_1}{\la_2}$ large.
\item[\rm b)] If $\be\le\tau_0\mu_2$ then \eqref{eq:system}with $N=3$ has a unique positive solution for $\frac{\la_1}{\la_2}>0$ small. If $\be\le\tau_0\mu_1$ then \eqref{eq:system} with $N=3$ has a unique positive solution for $\frac{\la_1}{\la_2}$ large.
\end{itemize}
\et

Theorems~\ref{thm:nonexist-2} and \ref{thm:sirakov} will be proved in Section~\ref{sec:asympt}.

\br\lab{rem:unique}
It is known and easy to see (cf.\ \cite{BartschWang2006, Sirakov2007}) that the problem
\beq\lab{eq:basic}
\begin{cases}
-\De u+u = \mu_1 u^3+\be uv^2,\quad &\hbox{in}\;\R^3,\\
-\De v+v = \mu_2 v^3+\be u^2v,\quad&\hbox{in}\;\R^3,\\
u,v>0,&\hbox{in}\;\R^3.
\end{cases}
\eeq
has no solution in the regime $\be\in [\min\{\mu_1,\mu_2\},\max\{\mu_1,\mu_2\}]$, if $\mu_1\ne\mu_2$. On the other hand, for $\be\in (0,\min\{\mu_1,\mu_2\})\cup (\max\{\mu_1,\mu_2\},+\infty)$ it is also easy to see that
\[
  u_\be(x) = \sqrt{\frac{\be-\mu_2}{\be^2-\mu_1\mu_2}}U(x),\quad
  v_\be(x) = \sqrt{\frac{\be-\mu_1}{\be^2-\mu_1\mu_2}}U(x)
\]
solve \eqref{eq:basic}. The solution $(u_\be,v_\be)$ is nondegenerate in the space $H_{rad}^1(\R^3,\R^2)$; see \cite[Lemma 2.2]{DancerWei2009}. Sirakov \cite[Remark~2]{Sirakov2007}) conjectured that, up to translations, $(u_\be,v_\be)$ is the unique positive solution of \eqref{eq:basic}. Wei and Yao \cite[Theorem~4.1, Theorem~4.2]{Wei2012} proved this conjecture for $\be>\max\{\mu_1,\mu_2\}$ and for $0<\be<\be_0$ close to $0$. Chen and Zou \cite[Theorem~1.1]{Chen2013a} proved the conjecture in case $\be'_0<\be<\min\{\mu_1,\mu_2\}$ close to $\min\{\mu_1,\mu_2\}$. The remaining range $\be\in [\be_0,\be'_0]$ is open up to now.
\er

%\br\lab{20180720-r1}
%Since $\be_1$ given by \eqref{20180720-e1} depends on $a,b$ (indeed depends on $\frac{a}{b}$), denote it by $\be_{1}=\be_{1}(a,b)$, then we also note that
%$$\sup_{\frac{a}{b}>0}\be_1(a,b)<\sqrt{\mu_1\mu_2}.$$
%To show this, we firstly  see that $\displaystyle\sup_{\frac{a}{b}>0}\be_1(a,b)$ is attained since $\be_1(a,b)\to 0$ as $\frac{a}{b}\to 0$ or $\frac{a}{b}\to \infty$.
%If the inequality above does not hold, then there exist some $a>0,b>0$ such that $\be_1(a,b)\geq \sqrt{\mu_1\mu_2}$, then the monotonicity implies that
%$$\max\left\{\frac{1}{a^2\mu_1^2},\frac{1}{b^2\mu_2^2}\right\}\leq \frac{1}{a^2(\mu_1+\sqrt{\mu_1\mu_2})^2}+\frac{1}{b^2(\mu_2+\sqrt{\mu_1\mu_2})^2}.$$
%Without loss of generality, we may assume that $a^2\mu_1^2\geq b^2\mu_2^2$, then
%$$\frac{1}{a^2(\mu_1+\sqrt{\mu_1\mu_2})^2}+\frac{1}{b^2(\mu_2+\sqrt{\mu_1\mu_2})^2}\geq \frac{1}{b^2\mu_2^2},$$
%i.e.,
%$$\frac{b^2}{a^2(\mu_1+\sqrt{\mu_1\mu_2})^2}+\frac{1}{(\mu_2+\sqrt{\mu_1\mu_2})^2}\geq \frac{1}{\mu_2^2}.$$
%Inserting $\frac{b^2}{a^2}\leq \frac{\mu_1^2}{\mu_2^2}$, we have
%$$\frac{\mu_1^2}{\mu_2^2(\mu_1+\sqrt{\mu_1\mu_2})^2}+\frac{1}{(\mu_2+\sqrt{\mu_1\mu_2})^2}\geq \frac{1}{\mu_2^2},$$
%which implies
%$$\frac{\mu_1+\mu_2}{(\sqrt{\mu_1}+\sqrt{\mu_2})^2}\geq 1,$$
%a contradiction to $\mu_1>0,\mu_2>0$.
%%From this view point, our result in present will cover all $0<\be<\sqrt{\mu_1\mu_2}$ is better than \cite[Theorem 1.1]{Bartsch2016}.
%\er

%%%%%%%%%%%%%
\s{Some  Preliminaries}\lab{sec:prelim}
%%%%%%%%%%%%%
In this section we collect results that hold for more general $N$, not only for $N=3$. We write $|u|_p$ for the $L^p$-norm. Let us first recall two results from \cite{Bartsch2017}.

\bl\lab{20180712-l1}
Let $(u,v)$ be a solution to
\beq\lab{20180712-e1}
\begin{cases}
-\De u+\la_1u=\mu_1 u^3+\be uv^2\quad &\text{in $\R^N$,}\\
-\De v+\la_2v=\mu_2 v^3+\be u^2v\quad&\text{in $\R^N$,}\\
u\geq 0,v\geq 0&\text{in $\R^N$,}
\end{cases}
\eeq
with $N\le3$. If $\la_1>0$ then  there exists $\al,\ga>0$ such that
$$u(x)\leq \al e^{-\sqrt{1+\ga|x|^2}}\;\hbox{for every }\;x\in \R^N.$$
%Similarly, if $\la_2>0$, we can obtain the similar decay property for $v$.
\el

Although only the case $N=3$ has been considered in \cite[Lemma 3.11]{Bartsch2017} the proof works verbatim for $N\le3$. The second result \cite[Lemma 3.12]{Bartsch2017} is a Liouville-type theorem.

\bl\lab{20180712-l2}
If $0\leq u\in H^1(\R^N)$ satisfies
$$-\De u+c(x)u\geq 0\;\hbox{in}\;\R^N,N\leq 3,$$
with $0\leq c(x)\leq Ce^{-C|x|}$ for some $C>0$, then $u\equiv 0$.
\el

\bp
The proof in \cite[Lemma 3.12]{Bartsch2017} for $N=3$ can be modified to cover $N\le2$ as follows. Suppose by contradiction that $u\not\equiv0$, hence $u>0$ by the strong maximum principle. Setting $v(x):= |x|^{-\al}$ for some $0 < \al \le \frac12$ there holds
\begin{align*}
 -\De v+c(x)v  &= \al(-\al+N-2)|x|^{-\al-2}+c(x)v\\
                        &\leq \al(-\al+N-2)|x|^{-\al-2}+C e^{-C|x|}|x|^{-\al}<0
\end{align*}
for every $|x|>r_0$ with $r_0$ large enough. Since $u>0$ in $\R^N$, there exists $C_0>0$ such that $u(x) \ge C_0 r_{0}^{-\al}$ for $|x|=r_0$. Now the comparison principle implies $u>C_0|x|^{-\al}$ in $\R^N\setminus B_{r_0}(0)$, hence $|u|_2 = \infty$, contradicting $u\in H^1(\R^N)$.
\ep

\bl\lab{lem:sol-properties}
Assume that $u,v\in H^1(\R^3)$ are positive and solve \eqref{eq:system} with $\mu_1,\mu_2>0$ and $\be\neq 0$. If in addition
$$
   \int_{\R^N}\big(\mu_1 u^4+\mu_2 v^4+2\be u^2v^2\big) > 0
$$
then $\la_1,\la_2>0$. Moreover, $u,v$ are radial functions (up to translation) and strictly radially decreasing if $\be>0$.
\el

\bp
We first observe that
$$
  |\nabla u|_2^2+\la_1|u|_2^2=\mu_1|u|_4^4+\be|uv|_2^2,\quad
  |\nabla v|_2^2+\la_2|v|_2^2=\mu_2|v|_4^4+\be|uv|_2^2,
$$
hence
\[%\beq\lab{20170313-e2}
  |\nabla u|_2^2+|\nabla v|_2^2
    = -(\la_1 |u|_2^2+\la_2|v|_2^2)+\big(\mu_1|u|_4^4+\mu_2|v|_4^4+2\be|uv|_2^2\big).
\]%\eeq
Now the Pohozaev identity
\[
\begin{aligned}
   &(N-2)\big(|\nabla u|_2^2+|\nabla v|_2^2\big)\\
   &\hspace{1cm}
     =-N\big(\la_1 |u|_2^2+\la_2 |v|_2^2\big)+\frac{N}{2}\big(\mu_1 |u|_4^4+\mu_2 |v|_4^4+2\be |uv|_2^2\big)
\end{aligned}
\]
implies
\[%\beq\lab{20170313-e4}
  \big(\la_1 |u|_2^2+\la_2 |v|_2^2\big)=\frac{4-N}{4}\big(\mu_1 |u|_4^4+\mu_2 |v|_4^4+2\be |uv|_2^2\big)>0.
\]%\eeq
Therefore without loss of generality we may assume $\la_1>0$. Then $u(x)$ decays exponentially at infinity according to Lemma~\ref{20180712-l1}. If $\la_2\leq 0$ we distinguish by the sign of $\be$. In the case $\be<0$, we have
$$
   -\De v+(-\be u^2)v = \mu_2 v^3-\la_2v \geq 0.
$$
Then $0\leq c(x):=-\be u^2\leq C e^{-C|x|}$ and $-\De v+c(x)v\geq 0$, hence $v\equiv 0$ by Lemma \ref{20180712-l2}. In the case $\be\geq0$, we have
$$
   -\De v \geq \mu_2 v^3\ \text{ in $\R^N$ and $v\geq 0$.}
$$
Now the classical Liouville-type theorem from \cite{GidasSpruck.1981} yields $v\equiv 0$, a contradiction. The last statement is due to \cite[Theorem~1]{BuscaSirakov2001}.
\ep

Let $S$ be the sharp constant for the embedding $H^1(\R^N)\hookrightarrow L^4(\R^N)$, i.e.\
\beq\lab{20170313-xe6}
  S|u|_{4}^{2} \le \big(|\nabla u|_2^2+|u|_2^2\big)\quad\text{for all $u\in H^1(\R^N)$,}
\eeq
and
\beq\lab{20170313-xe7}
  S=\left(|\nabla U|_2^2+|U|_2^2\right)^{\frac{1}{2}}=|U|_4^2
\eeq
where $U$ is the positive radial solution of \eqref{eq:def-U}. As in \cite[(1.6)]{BartschWangWei2007} we introduce the function $\tau:\R^+\to\R^+$ defined by
\beq\lab{eq:tau}
  \tau(s) := \inf_{\phi\in H^1(\R^N)\setminus\{0\}} \frac{\int_{\R^N} \left(|\nabla\phi|^2+s\phi^2\right)}{\int_{\R^N} U^2\phi^2}.
\eeq

\bl\lab{lem:tau}
\begin{itemize}
\item[\rm a)] The infimum $\tau_0$ in \eqref{eq:def-tau0} and the infimum in \eqref{eq:tau} are achieved by unique positive radial functions (and their scalar multiples).
\item[\rm b)] $\tau\in\cc^0(\R^+,\R^+)$ is strictly increasing and satisfies: $\tau(1)=1$, $\tau(s)\to\tau_0$ as $s\to0$, $\tau(s)\to\infty$ as $s\to\infty$.
\end{itemize}
\el

\bp
a) follows in a standard way from the compactness of the embedding $\cd_{0,rad}^{1,2}\hookrightarrow L^2(U^2dx)$ and symmetrization. The positive radial minimizer $\phi_s$, $s\ge0$, is the first eigenfunction of the eigenvalue problem $-\De\phi + s\phi = \la U^2\phi$. We choose $\phi_s$ to be normalized in $L^2(U^2dx)$.

b) We have for $s_1>s_2>0$:
$$
   \tau(s_2) < |\nabla \phi_{s_1}|_2^2+s_2|\phi_{s_1}|_2^2 < |\nabla \phi_{s_1}|_2^2+s_1|\phi_{s_1}|_2^2 = \tau(s_1),
$$
hence $\tau(s)$ is strictly increasing.

In order to prove the continuity consider a sequence $s_n\to s>0$. Clearly the minimizers $\phi_{s_n}$ are bounded, hence up to a subsequence $\phi_{s_n}\weakto \phi$ in $H^1(\R^N)$, and $\phi_{s_n}\to \phi$ in $L^2(U^2dx)$. This implies:
\begin{align*}
\tau(s) &\leq |\nabla \phi|_2^2+s|\phi|_2^2
     \leq \liminf_{n\to\infty}\left(|\nabla \phi_{s_n}|_2^2+s|\phi_{s_n}|_2^2\right)
     = \liminf_{n\to\infty}\tau(s_n) \\
   &\le \limsup_{n\to \infty}\tau(s_n)
     \leq \limsup_{n\to\infty} |\nabla \phi_s|_2^2+s_n|\phi_s|_2^2
     = |\nabla \phi_s|_2^2+s|\phi_s|_2^2 = \tau(s)
\end{align*}
Thus, $\tau(s_n) \to \tau(s)$ and $\phi=\phi_s$, so $\tau$ is continuous. Moreover, for $s>0$ we have $\phi_{s_n}\to\phi_s$ in $H^1(\R^N)$ because
$$
   |\nabla \phi_{s_n}|_2^2+s|\phi_{s_n}|_2^2 = \tau(s_n)+o(1) \to \tau(s) = |\nabla \phi_s|_2^2+s|\phi_s|_2^2.
$$

The identity $\tau(1)=1$ is obvious because by definition $U>0$ is an eigenfunction of $-\De\phi + \phi = \la U^2\phi$ associated to the eigenvalue $\la=1$.

Next  we observe that $\int_{\R^N} U^2\phi_s^2dx = 1$ and $U\in L^\infty(\R^N)$ imply $|\phi_s|_2\ge \ka > 0$ uniformly in $s$, hence
$$\tau(s) = |\nabla \phi_s|_2^2+s|\phi_s|_2^2 \ge s\ka^2 \to \infty\quad\text{as $s\to\infty$.}$$

In order to prove $\tau(s)\to\tau_0$ as $s\to0$ assume to the contrary that there exists $\de>0$ so that
$$
   \tau(s)\ge \tau_0+\de,\quad\text{for all $s>0$.}
$$
We choose a smooth cut-off function $\chi:\R\to[0,1]$ that is decreasing and satisfies
\[
  \chi(r) = \bcs 1 &\text{if $r\le 1$;}\\ 0 &\text{if $r\geq 2$.}\ecs
\]
Setting $\chi_R:\R^N\to\R$, $\chi_R(x) = \chi(|x|/R)$ we have for $R>0$ large that
$$
  \frac{|\nabla (\phi_0\chi_R)|_2^2}{\int_{\R^N}U^2(\phi_0\chi_R)^2 dx} < \tau_0+\frac12\de.
$$
This implies for $s$ close to $0$ the contradiction:
$$
   \tau(s) \le \frac{|\nabla(\phi_0\chi_R)|_2^2+s|\phi_0\chi_R|_2^2}{\int_{\R^N}U^2(\psi_0\chi_R)^2 dx}
    < \tau_0+\de
$$
\ep

%%%%%%%%%%%%%%%%%%%%%%%
\s{Global branches of solutions}\label{sec:global}
%%%%%%%%%%%%%%%%%%%%%%%
We consider a special case of \eqref{eq:system} , namely
\beq\lab{eq:bif-problem}
\begin{cases}
-\De u+\la u=\mu_1 u^3+\be v^2u\;\quad&\hbox{in}\;\R^3,\\
-\De v+v=\mu_2v^3+\be u^2v&\hbox{in}\;\R^3.
\end{cases}
\eeq
A straightforward computation shows the relation to \eqref{eq:main}.

\bl\lab{corres-1}
If $(u_\la,v_\la)$ is a solution of \eqref{eq:bif-problem} with
\beq\lab{eq:compare-norms}
  \frac{|u_\la|_2}{a} = \frac{|v_\la|_2}{b} =: \al
\eeq
then
\[
  u(x) = \al^2u_\la(\al^2x)\quad\text{and}\quad
  v(x) = \al^2 v_\la(\al^2x)
\]
solve \eqref{eq:main} with $\la_1=\la\al^4$ and $\la_2=\al^4$.
\el

\br
Clearly the converse holds in Lemma~\ref{corres-1}. If $(u,v)$ solves \eqref{eq:main} then
\[
  u_\la(x) = \sqrt{\la_2}u(\sqrt{\la_2}x)\quad\text{and}\quad
  v_\la(x) = \sqrt{\la_2}v(\sqrt{\la_2}x)
\]
solve \eqref{eq:bif-problem} with $\la=\frac{\la_1}{\la_2}$ and such that \eqref{eq:compare-norms} holds.
\er

Recall the solution $U$ of \eqref{eq:def-U}. Setting
\[
  U_{\la,\mu}(x) = \frac{\sqrt{\la}}{\sqrt{\mu}} U(\sqrt{\la}x)
\]
one easily checks that $(U_{\la,\mu_1},0)$ and $(0,U_{1,\mu_2})$ solve \eqref{eq:bif-problem}. These are called semitrivial solutions in the literature. We fix $\mu_1,\mu_2>0$ and consider $\la$ and $\be$ as parameters in \eqref{eq:bif-problem}. Then we have two families of semitrivial solutions of \eqref{eq:bif-problem}:
\[
  \ct_1 = \{(\la,\be,U_{\la,\mu_1},0):\ \la,\be>0\} \quad\text{and}\quad
  \ct_2 = \{(\la,\be,0,U_{1,\mu_2}):\ \la,\be>0\}.
\]
Clearly we also have the family $\ct_0:=\{(\la,\be,0,0):\ \la,\be>0\}$ of trivial solutions. Setting $E=H_{rad}^1(\R^3,\R^2)$ and $\P=\{(u,v)\in E: u,v\ge0\}$ for the positive cone, there holds $\ct_1,\ct_2 \subset X:= (\R^+)^2\times\P$; here $\R^+=(0,\infty)$. %Given $\be>0$ we write $X^\be := \R^+\times\{\be\}\times\P$ and use the notation $M^\be:=M\cap X^\be$ for subsets $M\subset X$.

We are interested in the set
\[
  \cs = \{(\la,\be,u,v)\in X: (\la,\be,u,v) \text{ solves }\eqref{eq:bif-problem},\ u,v>0\}
\]
of nontrivial positive solutions. Let us introduce the function
\beq\lab{eq:rho}
  \rho: \cs\to\R^+,\quad (\la,\be,u.v)\mapsto \frac{|u|_2}{|v|_2}.
\eeq
Lemma~\ref{corres-1} implies the following corollary which is the basic tool of our approach to finding normalized solutions.

\bc\lab{cor:corres-1}
If $\frac{a}{b}\in\rho(\cs^\be)$ then \eqref{eq:main} has a solution.
\ec

For the proof of Theorem~\ref{Main-th} it remains to get information about the image $\rho(\cs^\be)$. We shall approach this using continuation methods and bifurcation theory. First we investigate the solutions bifurcating from $\ct_1$ and $\ct_2$. Since we are interested in global bifurcation we reformulate \eqref{eq:bif-problem}. For $\la,\be>0$ we define a map $\A_{\la,\be}:\P\to\P$ by
\[
  \A_{\la,\be}(u,v) := \left((-\De+\la)^{-1}(\mu_1u^3+\be v^2u),(-\De+1)^{-1}(\mu_2v^3+\be u^2v)\right).
\]
As a consequence of the compact embedding $H_{rad}^1(\R^3)\hookrightarrow L^4(\R^3)$ the map
\[
  \A: X \to \P,\ \A(\la,\be,u,v) = \A_{\la,\be}(u,v),
\]
is completely continuous. Clearly fixed points of $\A_{\la,\be}$ correspond to solutions of \eqref{eq:bif-problem}. The set of bifurcation points can be explicitly determined. In order to describe it we define the functions
\beq\lab{eq:beta_i}
  \be_1(\la) = \mu_1\tau(1/\la) \quad\text{and}\quad \be_2(\la) = \mu_2\tau(\la) \quad\text{for $\la>0$}
\eeq
with $\tau$ from \eqref{eq:tau}. Using the fixed point index in the cone $\P$, denoted by $\ind_{\P}$, the following results have been proved in \cite{BartschWangWei2007}.

\bo\lab{prop:bif}
\begin{itemize}
\item[\rm a)] The map $\cs\to\R^+\times\R^+$, $(\la,\be,u,v)\mapsto (\la,\be)$ is proper, i.e.\ inverse images of compact sets are compact.

\item[\rm b)] $\overline\cs\cap\ct_1 = \big\{(\la,\be,U_{\la,\mu_1},0):\ \la>0,\ \be=\be_1(\la)\big\} =: \cb_1$

\item[\rm c)] $\overline\cs\cap\ct_2 = \big\{(\la,\be,0,U_{1,\mu_2}):\ \la>0,\ \be=\be_2(\la)\big\} =: \cb_2$

\item[\rm d)] For $\la,\be>0$ fixed we have
\[
  \ind_{\P}\big(\A_{\la,\be},(U_{\la,\mu_1},0)\big) = \bcs -1& \be<\be_1(\la)\\ 0 &  \be>\be_1(\la) \ecs
\]
and
\[
  \ind_{\P}\big(\A_{\la,\be},(0,U_{1,\mu_2})\big) = \bcs -1& \be<\be_2(\la)\\ 0 &  \be>\be_2(\la) \ecs
\]
\end{itemize}
\eo

In fact, in \cite{BartschWangWei2007} problem \eqref{eq:system} has been treated as a 5-parameter problem with parameters $(\la_1,\la_2,\mu_1,\mu_2,\be)\in(\R^+)^5$. The statement in \cite[Theorem~1.1]{BartschWangWei2007} about which part of $(\R^+)^5$ is covered by $\cs$ is not correct.

As a consequence of Proposition~\ref{prop:bif} there exist global two-dimensional continua $\cs_i\subset\cs$ bifurcating from $\ct_i$ so that $\overline\cs_i\cap\ct_i=\cb_i$, $i=1,2$. Using the analyticity of $\A$ it can be proved that $\cs$ and $\cs_i$ are two-dimensional manifolds except for one-dimensional subsets where secondary bifurcation takes place, but we do not need this. The global property of $\cs_i$ can be formulated as in \cite{AlexanderAntman1981}. This is somewhat technical and not needed here because we are interested in the case of prescribed $\be>0$. We will only use the standard Rabinowitz alternative for one-parameter global bifurcation.

As a corollary of Lemma~\ref{lem:tau} we obtain the following properties of the functions $\be_i$ defined in \eqref{eq:beta_i}.

\bc\label{cor:beta_i}
\begin{itemize}
\item[\rm a)] The function $\be_1$ is strictly decreasing and $\be_2$ is strictly increasing in $\la\in\R^+$.
\item[\rm b)] $\be_1(\la) \to \bcs \infty & \la\to0\\ \mu_1\tau_0 & \la\to\infty \ecs$
\item[\rm c)] $\be_2(\la) \to \bcs \mu_2\tau_0 & \la\to0\\ \infty & \la\to\infty \ecs$
\item[\rm d)] There exists a unique $\la^*>0$ such that $\be_1(\la^*) = \be_2(\la^*) =: \be^*$.
\end{itemize}
\ec
$$
   \includegraphics[height=3 in]{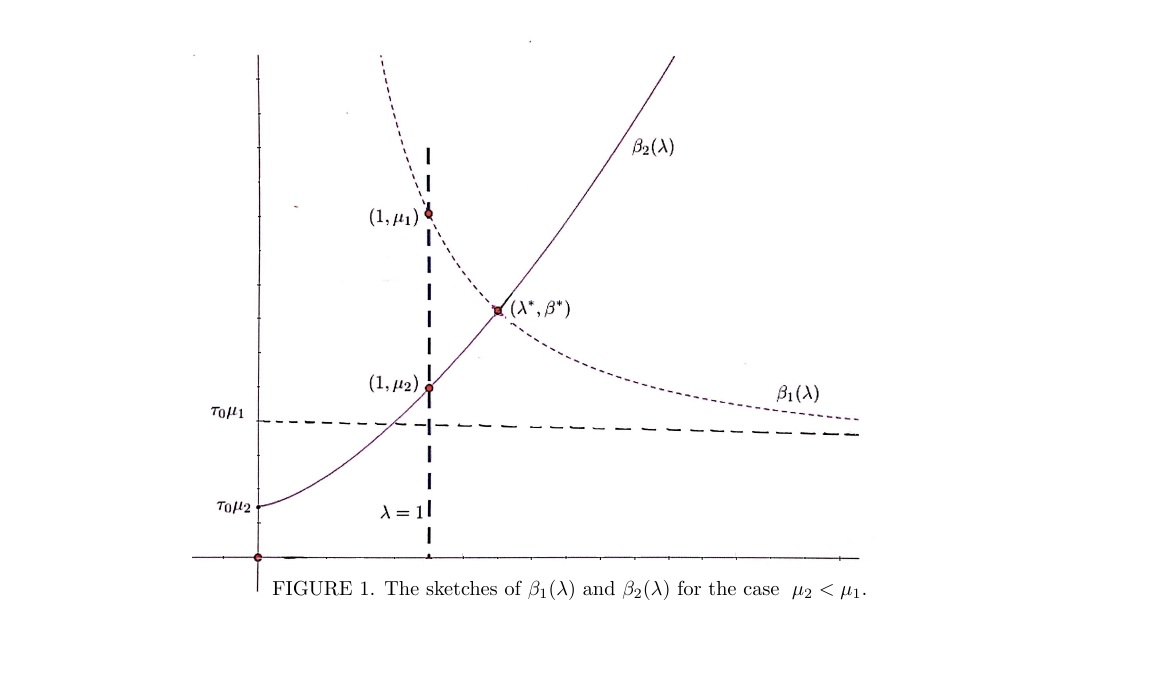}
$$

Now we deduce the global properties of the solutions bifurcating from $\ct_i$ that we need for $\be>0$ fixed. We set $\ell_i=\be_i^{-1}:(\mu_i\tau_0,\infty) \to \R^+$ for $i=1,2$, define $X^\be := \R^+\times\{\be\}\times\P$ for $\be>0$, and write $P_1: X \to\R^+$ for the projection onto the $\la$-component. For subsets $M\subset X$ we use the notation $M^\be:=M\cap X^\be$. The closure $\overline{M}$ of $M\subset X$ has to be understood in the relative topology of $X$.

\bo\lab{prop:global}
\begin{itemize}
\item[\rm a)] There is no bifurcation from the set $\ct_0=(\R^+)^2\times\{(0,0)\}$ of trivial solutions, i.e.\ $\overline{\cs}\cap\ct_0=\emptyset$.
\item[\rm b)] If $\be \le \tau_0\min\{\mu_1,\mu_2\}$ then $\overline{\cs^\be}\cap\ct_i^\be=\emptyset$, $i=1,2$.
\item[\rm c)] If $\mu_1\tau_0 < \be \le \mu_2\tau_0$ then there exists a connected component $\cs_1^\be \subset \cs^\be$ with $\overline{\cs_1^\be}\cap\ct_1^\be = \{(\ell_1(\be),\be,U_{\la,\mu_1},0)\}$. The projection $P_1(\cs_1^\be)$ contains the interval $(0,\ell_1(\be))$ or the interval $(\ell_1(\be),\infty)$. There is no bifurcation from $\ct_2^\be$ in $X^\be$.
\item[\rm d)]  If $\mu_2\tau_0 < \be \le \mu_1\tau_0$ then there exists a connected component $\cs_2^\be \subset \cs^\be$ with $\overline{\cs_2^\be}\cap\ct_2^\be = \{(\ell_2(\be),\be,0,U_{1,\mu_2})\}$. The projection $P_1(\cs_2^\be)$ contains the interval $(0,\ell_2(\be))$ or the interval $(\ell_2(\be),\infty)$. There is no bifurcation from $\ct_1^\be$ in $X^\be$.
\item[\rm e)] If $\be > \tau_0\max\{\mu_1,\mu_2\}$ then there exist connected sets $\cs_i^\be \subset \cs^\be$, $i=1,2$, with $\overline{\cs_1^\be}\cap\ct_1^\be = \{(\ell_1(\be),\be,U_{\la,\mu_1},0)\}$ and $\overline{\cs_2^\be}\cap\ct_2^\be = \{(\ell_2(\be),\be,0,U_{1,\mu_2})\}$. If $\cs_1^\be\cap\cs_2^\be\ne\emptyset$ then $\cs_1^\be=\cs_2^\be$. If this is not the case then $P_1(\cs_1^\be)$ contains the interval $(0,\ell_1(\be))$ or the interval $(\ell_1(\be),\infty)$, and $P_1(\cs_2^\be)$ contains the interval $(0,\ell_2(\be))$ or the interval $(\ell_2(\be),\infty)$.
\end{itemize}
\eo

\bp
a) This is clear since $(0,0)$ is a nondegenerate solution of \eqref{eq:bif-problem} for all $(\la,\be)\in(\R^+)^2$.

b) As a consequence of Corollary~\ref{cor:beta_i} there is no $\la>0$ with $\beta_1(\la)=\be$ or $\be_2(\la)=\be$.

c) Here Corollary~\ref{cor:beta_i} implies that there exists $\la_1=\ell_1(\be)>0$ with $\beta_1(\la_1)=\be$ but there is no $\la_2>0$ with $\beta_2(\la_2)=\be$. Therefore there exists a connected set $\cs_1^\be \subset \big((id-\A)^{-1}(0)\cap X^\be\big)\setminus\ct_1$ with $\overline{\cs_1^\be}\cap\ct_1^\be = \{(\ell_1(\be),\be,U_{\la,\mu_1},0)\}$ and which satisfies the classical Rabinowitz alternative. It cannot return to $\ct_1^\be$ because there is no second bifurcation point on $\ct_1^\be$. Therefore it must be unbounded. Since there is no bifurcation from $\ct_0$ and $\ct_2$ we deduce that $\overline{\cs_1^\be}\cap\ct_i^\be=\emptyset$, $i=0,2$, hence $\cs_1^\be\subset\cs$. Now Proposition~\ref{prop:bif}~a) implies that the only way for $\cs_1^\be$ to be unbounded is that $P_1(\cs_1^\be)$ contains the interval $(0,\ell_1(\be))$ or the interval $(\ell_1(\be),\infty)$. To be careful, if $P_1(\cs_1^\be)$ contains the interval $(0,\ell_1(\be))$ then $\cs_1^\be$ is already unbounded in the sense of the Rabinowitz alternative because we only consider the parameter range $\la\in\R^+$. It is not necessary that the $(u,v)$-component becomes unbounded in $\cs_1^\be$.

d) The proof is analogous to the one of c).

e) As in the proof of c) and d) there exist connected sets $\widetilde\cs_i^\be\subset  \big((id-\A)^{-1}(0)\cap X^\be\big)\setminus\ct_i$ bifurcating from $\ct_i$ which satisfy the Rabinowitz alternative. If the closure of $\widetilde\cs_1^\be$ intersects $\ct_2^\be$ then $\widetilde\cs_1^\be$ contains $\ct_2$ and the connected set of nontrivial solutions bifurcating from $\ct_2$. This implies that
\[
  \cs_1^\be := \widetilde\cs_1^\be \cap \cs = \widetilde\cs_1^\be \setminus \ct_2^\be
                  = \widetilde\cs_2^\be \setminus \ct_1^\be= \widetilde\cs_2^\be \cap \cs  =:  \cs_2^\be
\]
connects $\ct_1^\be$ and$\ct_2^\be$. Analogously this holds if the closure of $\widetilde\cs_2^\be$ intersects $\ct_1^\be$.

It remains to consider the case where the closure of $\widetilde\cs_i^\be$ does not intersect $\ct_{3-i}^\be$ for $i=1,2$. Then $\cs_i^\be := \widetilde\cs_i^\be \subset \cs^\be$ is unbounded in the sense of c) and d), i.e.\ $P_1(\cs_i^\be)$ contains the interval $(0,\ell_i(\be))$ or the interval $(\ell_i(\be),\infty)$, $i=1,2$.
\ep

\br
Using analytic bifurcation theory one can prove that the sets $\cs_i^\be$ are smooth curves except for a discrete subset of singular points. One can also apply the Crandall-Rabinowitz theorem about bifurcation from simple eigenvalues to see that $\cs_i^\be$ is a curve near the bifurcation point. These results are not needed here.
\er

As a corollary we obtain a first major building block of the proof of Theorem~\ref{Main-th}.

\bc\lab{cor:cs1=cs2}
If $\be > \max\{\mu_1\tau_0,\mu_2\tau_0\}$ and $\cs_1^\be\cap\cs_2^\be\ne\emptyset$ then problem~\eqref{eq:main} has a solution for every $a,b>0$.
\ec

\bp
Recall the function $\rho$ from \eqref{eq:rho}. By definition there exist $(\la_n,\be,u_n,v_n) \in \cs_1^\be$ such that $(\la_n,\be,u_n,v_n) \to (\ell_1(\be),\be,U_{\ell_1(\be),\mu_1},0)\}$, hence $\rho(\la_n,\be,u_n,v_n) \to \infty$ as $n\to\infty$. And as a consequence of Proposition~\ref{prop:global}~e) there exist $(\la'_n,\be,u'_n,v'_n) \in \cs_1^\be$ such that $(\la'_n,\be,u'_n,v'_n) \to (\ell_2(\be),\be,0,U_{1,\mu_2})$, hence $\rho(\la'_n,\be,u'_n,v'_n) \to 0$ as $n\to\infty$. Since $\cs_1^\be$ is connected it follows that $\rho$ is onto. Now the result follows from Corollary~\ref{cor:corres-1}.
\ep

In addition to the global continua bifurcating from $\ct_1$ and $\ct_2$ there exists a third global continuum $\widetilde\cs\subset\cs$. In order to see this recall that for $\la=1$ and $\be\in(0,\be_0)$ close to $0$ the problem \eqref{eq:bif-problem} has precisely four solutions in $\P$: the trivial solution $(0,0)$, the semitrivial solutions $(U_{1,\mu_1},0)$, $(0,U_{1,\mu_2})$, and a unique nontrivial solution $(u_\be,v_\be)$ which satisfies $(u_\be,v_\be)\to(U_{1,\mu_1},U_{1,\mu_2})$ as $\be\to0$; see Remark~\ref{rem:unique}. The map
\[
  (0,\be_0)\to\P, \quad\be\mapsto (u_\be,v_\be),
\]
is smooth by the implicit function theorem applied at $(U_{1,\mu_1},U_{1,\mu_2})$.

\bo\lab{prop:cs_0}
For $\be\in(0,\be_0)$ there holds $\ind_\P(\A_{1,\be},(u_\be,v_\be))=1$.
\eo

\bp
The solution $(U_{1,\mu_1},U_{1,\mu_2})$ of \eqref{eq:bif-problem} with $\la=1$ and $\be=0$ has Morse index $2$ as critical point of $J$, with negative eigenspace spanned by $(U_{1,\mu_1},0),(0,U_{1,\mu_2})\in\P$. The Poincar\'e-Hopf theorem in convex sets \cite[Theorem~1.5]{BartschDancer2009} implies
\[
  \ind_\P(\A_{1,0},(U_{1,\mu_1},U_{1,\mu_2})=(-1)^2=1.
\]
Now the proposition follows from the homotopy invariance of the fixed point index.
\ep

The homotopy invariance of the fixed point index allows to continue the solutions $(u_\be,v_\be)$ to other parameter values in $(\R^+)^2$. We define $\widetilde\cs\subset\cs$ to be the connected component of $\cs$ containing the nontrivial solutions $(1,\be,u_\be,v_\be)$ for $\be>0$ small. As a corollary of Proposition~\ref{prop:cs_0} we obtain the following.

\bc\lab{cor:beta-small}
If $\be \le \tau_0\min\{\mu_1,\mu_2\}$ then there exists a connected set $\cs_0^\be \subset \cs^\be\cap\widetilde\cs$ such that $P_1(\cs_0^\be) = \R^+$.
\ec

\bp
Let $\co\subset X\setminus(\cs\cup\cb_1\cup\cb_2)$ be an open neighborhood of
\[
  \ct_0\cup (\ct_1\setminus\cb_1) \cup (\ct_2\setminus\cb_2) \subset X\setminus(\cs\cup\cb_1\cup\cb_2)
\]
such that $\cs\cap\overline{\co}=\emptyset$. For $\la,\be>0$ we set $\co_{\la,\be} := \{(u,v)\in\P:(\la,\be,u,v)\in\co\}$. By definition the nontrivial fixed points of $\A_{\la,\be}$ are contained in $\Om_{\la,\be} :=B_R(0) \setminus \overline{\co_{\la,\be}}$ for $R>R(\la,\be)$ large. This a bounded and open subset of $\P$. Proposition~\ref{prop:cs_0} and the homotopy invariance of the fixed point index imply for $\be \le \min\{\tau_0\mu_1,\tau_0\mu_2\}$ and $\be'\in(0,\be_0)$:
\[
  \ind_\P(\A_{\la,\be},\Om_{\la,\be}) = \ind_\P(\A_{\la,\be'},\Om_{\la,\be'})
	   = \ind_\P(\A_{1,\be'},\Om_{1,\be'}) = 1
\]
The result follows from the continuation principle.
\ep

Observe that $\cs_0^\be$ may differ from $\widetilde\cs^\be=\widetilde\cs\cap X^\be$ because the latter may not be connected.

We may also use Proposition~\ref{prop:cs_0} to compute the global fixed point index of all positive solutions of \eqref{eq:bif-problem}, for each $\la,\be>0$. Observe that according to Proposition~\ref{prop:bif}~a) for $\la,\be>0$ there exists $R(\la,\be)>0$ such that the positive solutions of \eqref{eq:bif-problem} are bounded by $R(\la,\be)$. Therefore the fixed point index
\[
  i_\infty(\la,\be) = \ind_\P(\A_{\la,\be},B_R(0))
\]
is well defined and independent of $R>R(\la,\be)$. Applying the homotopy invariance of the fixed point index and Proposition~\ref{prop:bif}~a) again, we also see that $i_\infty := i_\infty(\la,\be)$ is independent of $\la,\be>0$.

\bo\lab{prop:i_infty}
$i_\infty=0$
\eo

\bp
We compute $i_\infty(\la,\be)$ for $\la=1$ and $\be\in(0,\be_0)$. Then $i_\infty = i_\infty(1,\be)$ is the sum of the local indices at the four solutions $(0,0)$, $(U_{1,\mu_1},0)$, $(0,U_{1,\mu_2})$, $(u_\be,v_\be)$. From \cite[Theorem~1.5]{BartschDancer2009} it follows that
\[
  \ind_\P\big(\A_{1,0},(0,0)\big) = 1.
\]
Propositions~\ref{prop:bif} and \ref{prop:cs_0} imply for $\be\in(0,\be_0)$:
\[
\begin{aligned}
  i_\infty &= \ind_\P\big(\A_{1,\be},(0,0)\big) + \ind_\P\big(\A_{1,\be},(U_{1,\mu_1},0)\big)
	            + \ind_\P\big(\A_{1,\be},(0,U_{1,\mu_1})\big)\\
	     &\hspace{1cm} + \ind_\P\big(\A_{1,\be},(u_\be,v_\be)\big)
		= 1-1-1+1 = 0
\end{aligned}
\]
\ep

%%%%%%%%%%%%%%%%%%%%%%%%%%%%%%%%%%%%%%%%%
\s{Asymptotic behavior of positive solutions for $\la\to0$ or $\la\to\infty$}\lab{sec:asympt}
%%%%%%%%%%%%%%%%%%%%%%%%%%%%%%%%%%%%%%%%%
In this section we investigate the function
\[
  \rho:\cs\to\R^+,\quad \rho(\la,\be,u,v) = \frac{|u|_2}{|v|_2},
\]
from \eqref{eq:rho} as $\la\to0$ or $\la\to\infty$.

\bl\lab{lem:la-to-0}
Let $(u_{n},v_{n})$, $n\in\N$, be positive radial solutions to equation \eqref{eq:bif-problem} with $\la=\la_n\to0$. Then the following conclusions hold up to a subsequence.
\begin{itemize}
\item[\rm a)] $u_{n}(x)+v_{n}(x) \to 0$ as $|x|\to \infty$  uniformly in $n$.
\item[\rm b)]  $|u_{n}|_\infty\to 0$, $|v_{n}|_\infty\le C$, and $(u_{n},v_{n})\to(0,U_{1,\mu_2})$ in $\cc^2_{loc}(\R^N)\times \cc^2_{loc}(\R^N)$.
\item[\rm c)] $v_{n}\to U_{1,\mu_2}$ in $H^1(\R^N)$
\item[\rm d)] $|\nabla u_{n}|_2 = O(1)|u_{n}|_2$; if $u_{n}$ is unbounded in $H^1(\R^N)$, then $\rho(\la_n,\beta,u_{n},v_{n})\to \infty$.
\end{itemize}
\el

\bp
a) The proof in \cite[Step 2 in the proof of Theorem 1.1]{Chen2013a} is valid here.

b) A standard blow up argument as in \cite[Lemma 2.4]{DancerWei2009} shows that $|u_{n}|_\infty+|v_{n}|_\infty$ is bounded. If $\displaystyle\al := \liminf_{n\to \infty}u_{n}(0) > 0$ we consider
$$
   -\De \frac{u_{n}}{u_{n}(0)}+\la_n \frac{u_{n}}{u_{n}(0)}
    = \mu_1 u_{n}(0)^2 \left(\frac{u_{n}}{u_{n}(0)}\right)^3+\be v_{n}^{2} \frac{u_{n}}{u_{n}(0)}.
$$
Then $\frac{u_{n}}{u_{n}(0)}\to \tilde{u}$ as $n\to\infty$ along a subsequence, which is a nonnegative radial function satisfying
$$
   -\De \tilde{u}\ge \mu_1 \eps_0^2 \tilde{u}^3.
$$
Now \cite{GidasSpruck.1981} implies $\tilde{u}\equiv 0$, contradicting $\tilde{u}(0)=1$. Therefore $|u_{n}|_\infty\to 0$, hence $u_{n}\to 0$ in $C^2_{loc}(\R^N)$ along a subsequence. Since $v_{n}=(-\De +1)^{-1}(\mu_2 v_{n}^{3}+\be u_{n}^{2}v_{n})$ and $|u_{n}|_\infty\to 0$, we see that $|v_{n}|_\infty$ is bounded away from $0$. Then $\displaystyle\tilde{v} := \lim_{n\to \infty}v_{n}$ is a positive radial solution to
$$
   -\De v+v=\mu_2 v^3\,,\quad v(x) \to 0\text{ as $|x|\to\infty$,}
$$
which implies $\tilde{v}=U_{1,\mu_2}$ and $v_{n}\to U_{1,\mu_2}$ in $C^2_{loc}(\R^N)$.

c) It is standard to prove that $v_{n}(x)\to0$ exponentially and uniformly in $n$, so there exist $C,R>0$, independent of $n$ such that
$$
   v_{n}(x) \le C e^{-\frac{1}{2}|x|} \quad\text{for all $|x|>R$, all $n\in\N$.}
$$
As in b), or \cite[Step 3 in the proof of Theorem 1.1]{Chen2013a}, one sees that $v_{n}$ is bounded in $H^1(\R^N)$. Observe that this argument is not valid for $u_{n}$ because $\la_n\to 0$. Then we have, up to a subsequence:
$$
   \text{$v_{n}\weakto v$ in $H^1(\R^N)$, $v_{n}\to v$  in $L^4(\R^N)$, and $v_{n}\to v$ a.e. in $\R^N$,}
$$
which implies $v=U_{1,\mu_2}$. Now we recall that $|u_{n}|_\infty\to 0$, hence $\be|u_{n}v_{n}|_2^2\to 0$. Using
$$
   |\nabla v_{n}|_2^2+|v_{n}|_2^2=\mu_2|v_{n}|_4^4+\be|u_{n}v_{n}|_2^2
$$
and $v_{n}\to U_{1,\mu_2}$ in $L^4(\R^N)$, we deduce
$$
  |\nabla v_{n}|_2^2+|v_{n}|_2^2\to \mu_2|U_{1,\mu_2}|_4^4=|\nabla U_{1,\mu_2}|_2^2+|U_{1,\mu_2}|_2^2.
$$
This yields $v_{n}\to U_{1,\mu_2}$ in $H^1(\R^N)$.

d) Setting $|\nabla u_{n}|_2^2=\sigma_n |u_{n}|_2^2$ we have
$$(\sigma_n+\la_n)|u_{n}|_2^2=\mu_1|u_{n}|_4^4+\be |u_{n}v_{n}|_2^2.$$
Now a) and b) imply $\mu_1|u_{n}|_4^4+\be |u_{n}v_{n}|_2^2=O(1)|u_{n}|_2^2$, hence $|\nabla u_{n}|_2^2=O(1)|u_{n}|_2^2$. Thus if $u_{n}$ is unbounded in $H^1(\R^N)$ then $u_{n}$ must be unbounded in $L^2(\R^N)$ and
$\rho(\la_n,\beta,u_n,v_n) = \frac{|u_{n}|_2}{|v_{n}|_2}\to \infty$.
\ep

\bl\lab{lem:la-to-infty}
Let $(u_{n},v_{n})$, $n\in\N$, be positive radial solutions to equation \eqref{eq:bif-problem} with $\la=\la_n\to\infty$. Then $\bar{u}_{n}(x) := \frac{1}{\sqrt{\la_n}}v_{n}\left(x/\sqrt{\la_n}\right)$ and $\bar{v}_{n}(x) := \frac{1}{\sqrt{\la_n}}u_{n}\left(x/\sqrt{\la_n}\right)$ satisfy (along a subsequence):
\begin{itemize}
\item[\rm a)] $\bar{u}_{n}(x)+\bar{v}_{n}(x) \to 0$ as $|x|\to\infty$ uniformly in $n$.
\item[\rm b)]  $|\bar{u}_{n}|_\infty\to 0$, $|\bar{v}_{n}|_\infty\le C$, and $(\bar{u}_{n},\bar{v}_{n})\to(0,U_{1,\mu_1})$ in $\cc^2_{loc}(\R^N)\times \cc^2_{loc}(\R^N)$.
\item[\rm c)] $\bar{v}_{n}\to U_{1,\mu_1}$ in $H^1(\R^N)$
\item[\rm d)] $|\nabla\bar{u}_{n}|_2 = O(1)|\bar{u}_{n}|_2$; if $\bar{u}_{n}$ is unbounded in $H^1(\R^N)$ then $\rho(\la_n,\beta,u_{n},v_{n})\to \infty$.
\end{itemize}
\el

\bp
A direct computation shows that $(\bar{u}_{n},\bar{v}_{n})$ solve
$$
\bcs
   -\De u+\frac{1}{\la_n} u = \mu_2 u^3+\be uv^2 &\hbox{in}\;\R^N,\\
   -\De v+v = \mu_1v^3+\be vu^2 &\hbox{in}\;\R^N.
\ecs
$$
The result follows from Lemma \ref{lem:la-to-0} and
$$
   \rho(\la_n,\beta,u_{n},v_{n}) =\frac{|u_{n}|_2}{|v_{n}|_2}=\frac{|\bar{v}_{n}|_2}{|\bar{u}_{n}|_2}\to 0.
$$
\ep

Now we prove Theorems~\ref{thm:nonexist-2} and \ref{thm:sirakov}. Observe that $(u,v)$ is a positive solution to \eqref{eq:system} if and only if
$$
  \bar{u}(x) := \frac{1}{\sqrt{\la_2}}u\left(x/\sqrt{\la_2}\right), \quad
  \bar{v}(x) := \frac{1}{\sqrt{\la_2}}v\left(x/\sqrt{\la_2}\right),
$$
solve \eqref{eq:system} with $\la_1=\la$ and $\la_2=1$, i.e.\ \eqref{eq:bif-problem}. Therefore ist is sufficient to consider this case.\\

\begin{altproof}{Theorem~\ref{thm:nonexist-2}}
a) Arguing by contradiction suppose that for fixed $\be\ge\mu_2$ there exist a sequence $\la_n\to 0$ and positive solutions $(u_n,v_n)$ to \eqref{eq:bif-problem} with $\la=\la_n$. Then we have
\[
  \langle \nabla u_n,\nabla v_n\rangle +\la_n \int_{\R^N}u_nv_n=\mu_1\int_{\R^N}u_n^3v_n+\be\int_{\R^N}u_nv_n^3
\]
and
\[
  \langle \nabla u_n,\nabla v_n\rangle + \int_{\R^N}u_nv_n=\mu_2\int_{\R^N}v_n^3u_n+\be\int_{\R^N}v_nu_n^3.
\]
These identities yield
\[
  (1-\la_n)\langle \nabla u_n,\nabla v_n\rangle =\int_{\R^N}[(\be-\la_n\mu_2)v_n^3u_n+(\mu_1-\la_n\be)v_nu_n^3],
\]
which implies $\langle \nabla u_n,\nabla v_n\rangle >0$ for $n$ large enough. On the other hand, we also have
\[
  (1-\frac{\be}{\mu_2})\langle\nabla u_n,\nabla v_n\rangle + (\la_n-\frac{\be}{\mu_2})\int_{\R^N}u_nv_n
    =\int_{\R^N}(\mu_1-\frac{\be^2}{\mu_2})v_nu_n^3.
\]
Now $|u_n|_\infty\to 0$ by Lemma~\ref{lem:la-to-0}, so that
$$
   \int_{\R^N}(\mu_1-\frac{\be^2}{\mu_2})v_nu_n^3=o(1)\int_{\R^N}u_nv_n.
$$
In the case $\be=\mu_2$, we deduce
$$
   \frac{\be}{\mu_2}\int_{\R^N}u_nv_n=o(1)\int_{\R^N}u_nv_n,
$$
a contradiction. And if $\be>\mu_2$ we obtain
$$
  (1-\frac{\be}{\mu_2})\langle \nabla u_n,\nabla v_n\rangle=(\frac{\be}{\mu_2}+o(1))\int_{\R^N}u_nv_n>0,
$$
which implies $\langle \nabla u_n,\nabla v_n\rangle<0$ for $n$ large enough, a contradiction again.

b) This follows from a) using the transformation from the proof of Lemma \ref{lem:la-to-infty}.
\end{altproof}

Now we recall \cite[Lemma~2.3]{DancerWei2009}.

\bl\lab{lem:lin-problem}
The linearized problem
\[
\bcs
  \De \phi-\la\phi+3\mu_1u^2\phi+\be v^2\varphi+2\be uv\psi=0,\quad&x\in\R^N,\\
  \De\psi-\psi+3\mu_2v^2\psi+\be u^2\psi+2\be uv\phi=0,\quad&x\in\R^N,\\
  \varphi=\varphi(r), \phi=\phi(r),
\ecs
\]
has exactly a one-dimensional set of solutions for $\la>0$ and $\be= \be_{1}(\la)$, $(u,v)=(U_{\la,\mu_1},0)$ or $\be= \be_{2}(\la)$, $(u,v)=(0,U_{1,\mu_2})$.
\el

We have a similar result for $\la=0$.

\bl\lab{lem:lin-problem-la=0}
The linearized problem
\[
\bcs
  -\De \phi=\be U_{1,\mu_2}^2\phi,\quad&x\in\R^N,\\
  \De\psi-\psi+3\mu_2U_{1,\mu_2}^2\psi=0,\quad&x\in\R^N,\\
  \phi=\phi(r), \psi=\psi(r).
\ecs
\]
has only the zero solution if $0<\be \ne \tau_0 \mu_2$. If $\be=\tau_0\mu_2$ then the set of solutions has dimension one.
\el

\bp
It is well known that the eigenvalue problem
\[
  -\De\phi+\phi = \nu\mu_2\om_{1,\mu_2}^2\phi = \nu\om_{1,1}^2\phi
\]
has eigenvalues $\nu_1=1$, $\nu_2=\dots=\nu_{N+1}=3$,$\nu_k>3$ for $k\ge N+2$, and that the eigenfunctions corresponding to $\nu=3$ are not radial. It follows that $\psi=0$. If $\phi\not\equiv 0$ then $\phi>0$ by the maximum principle, and $\phi$ is a minimizer of $\be_{2}(0)=\mu_2\tau_0$. The result follows from Lemma~\ref{lem:tau}.
\ep

Now we return to study the asymptotic behavior of the positive solution for $\la$ small or large and improve on Lemmas~\ref{lem:la-to-0} and \ref{lem:la-to-infty}. And then give the proof of Theorem~\ref{thm:sirakov} to end this section.

\bl\lab{lem:la-to-0-or-infty}
a) Let $(u_{n},v_{n})$, $n\in\N$, be positive radial solutions of equation \eqref{eq:bif-problem} with $\la=\la_n\to0$. Then
$$
  \left(\frac{1}{\sqrt{\la_n}}u_n\left(x/\sqrt{\la_n}\right),v_n(x)\right) \to \big(U_{1,\mu_1}(x),U_{1,\mu_2}(x)\big)\quad
   \text{in $\cc_{loc}^2(\R^N)\times \cc_{loc}^2(\R^N)$.}
$$

b) Let $(u_{n},v_{n})$, $n\in\N$, be positive radial solutions of equation \eqref{eq:bif-problem} with $\la=\la_n\to\infty$. Then
$$
  \left(\frac{1}{\sqrt{\la_n}}u_n\left(x/\sqrt{\la_n}\right),v_n(x)\right) \to \big(U_{1,\mu_1}(x),U_{1,\mu_2}(x)\big)\quad
   \text{in $\cc_{loc}^2(\R^N)\times \cc_{loc}^2(\R^N)$.}
$$

\el

\bp
a) We first consider the case $\la_n\to0$.

{\sc Step 1:} $\liminf_{n\to\infty}\frac1{\sqrt{\la_n}}u_n(0)>0$.

We argue by contradiction and assume that $u_n(0)=o(1)\sqrt{\la_n}$, after passing to a subsequence. The function
$$
   \bar{u}_n(x):=\frac{1}{u_n(0)}u_n\left(x/\sqrt{\la_n}\right)
$$
solves
\beq\lab{eq:bar-u_n}
  -\De \bar{u}_n(x)+\bar{u}_n(x) = \frac{u_n(0)^2}{\la_n}\mu_1\bar{u}_n(x)^3+\be \bar{u}_n(x) \bar{v}_n(x)^2
\eeq
with
\[
   \bar{v}_n(x):= \frac{1}{\sqrt{\la_n}}v_n\left(x/\sqrt{\la_n}\right)\,.
\]
Observe that $\bar{u}_n\to \bar{u}$ in $\cc^0_{loc}(\R^N)$ along a subsequence and $\bar{u}(0)=1$ because $|\bar{u}_n|_\infty=\bar{u}_n(0)=1$. By Lemma~\ref{lem:la-to-0} we have $v_n\to U_{1,\mu_2}$ both in $H^1(\R^N)$ and in $C_{loc}^{2}$, and $v_n(x)\to0$ as $|x|\to\infty$ uniformly in $n$. It follows that $\bar{v}_n\to 0$ uniformly outside an arbitrary neighborhood of $0$. For a test function $h\in \cd(\R^N)$ and $\eps>0$, there exists $r_\eps$ such that
\[
%\begin{aligned}
  \int_{|x|\leq r_0} \big|\bar{u}_n\bar{v}_n^2(x) h(x)\big|dx
%    &\le |\tilde{v}_n|_3^2 \left(\int_{|x|\leq r_0}|h(x)|^3dx\right)^{\frac{1}{3}}\\
    \le |v_n|_3^2 \left(\int_{|x|\leq r_\eps}|h(x)|^3dx\right)^{\frac{1}{3}}<\frac{\eps}{2}.
%\end{aligned}
\]
Therefore  $\int_{\R^N} \bar{u}_n\bar{v}_n^2 h\,dx \to 0$. Testing \eqref{eq:bar-u_n} with $h$ we see that $\bar{u}_n\weakto 0$ in $H^1(\R^N)$, contradicting $\bar{u}_n\to \bar{u}$ in $\cc^0_{loc}(\R^N)$.

 {\sc Step 2:}  $\limsup_{n\to\infty}\frac1{\sqrt{\la_n}}u_n(0)<\infty$.

Assume by contradiction that $\sqrt{\la_n}=o(1)u_n(0)$, after passing to a subsequence. The function
$$
   \widetilde{u}_n(x)=\frac{1}{u_n(0)}u_n\big(\sqrt{\la_n}x/u_n(0)\big)
$$
satisfies $|\widetilde{u}_n|_\infty = \widetilde{u}_n(0) = 1$ and
$$
   -\De \widetilde{u}_n+\frac{\sqrt{\la_n}}{u_n(0)} \widetilde{u}_n \geq \mu_1 \widetilde{u}_n^3\quad\hbox{in $\R^N$.}
$$
Then $\tilde{u}_n\to \tilde{u}\geq 0$ in $C_{loc}^{2}(\R^N)$, along a subsequence, with $\tilde{u}(0)=1$, and $\tilde{u}$ satisfies
$$
   -\De \tilde{u} \geq \mu_1\tilde{u}^3\quad\hbox{in $\R^N$.}
$$
This implies $\tilde{u}\equiv 0$, a contradiction.

The conclusion about $v_n(x)$ has already been proved in Lemma~\ref{lem:la-to-0}.

{\sc Step 3:} $\displaystyle \bar{u}_n(x) := \frac{1}{\sqrt{\la_n}}u_n\left(x/\sqrt{\la_n}\right) \to U_{1,\mu_1}(x)$\ \ in $\cc_{loc}^2(\R^N)$

Observe that
%\beq%\lab{20190622-e1}
\[
\left\{
\begin{aligned}
-\De \bar{u}_n+\bar{u}_n
  &= \mu_1\bar{u}_n^3+\frac{\be}{\la_n} \bar{u}_nv_n^2\left(\,\cdot\,/\sqrt{\la_n}\right)
  &&\hbox{in}\;\R^N\\
-\De {v}_n+{v}_n
  &= \mu_2{v}_n^3+\be {v}_n\left(\sqrt{\la_n}\bar{u}_n\left(\sqrt{\la_n}\,\cdot\,\right)\right)^2 &&\hbox{in}\;\R^N.
\end{aligned}
\right.
\]
%\eeq
By  {\sc Step }1 and {\sc Step }2 we may assume that $\bar{u}_n\to \bar{u}\geq 0$ in $C_{loc}^2(\R^N)$ and $\bar{u}(0)>0$, hence $\bar{u}>0$ in $\R^N$. By $\la_n\rightarrow 0$, we may assume that $\la_n<1$ for all $n$. Recalling that there exist $C,R>0$, independent of $n$ such that
$$
  v_n(x)\leq C e^{-\frac{1}{2}|x|}\;\hbox{for all}\;|x|>R, \;\hbox{all}\;n\in \N,
$$
we have that
$$
  \frac{\be}{\la_n}v_n^2\left(x/\sqrt{\la_n}\right) \leq \be C^2 \frac{1}{\la_n}e^{-|x|/\sqrt{\la_n}}\;\hbox{ for all}\;|x|>R, \;\hbox{all}\;n\in \N.
$$
Fix $R>0$, then $\be C^2 \frac{1}{\la_n}e^{-R/\sqrt{\la_n}} \to 0$ as $n\rightarrow \infty$, which implies that
$$
   \frac{\be}{\la_n}v_n^2\left(x/\sqrt{\la_n}\right)<\frac{1}{2}\;\hbox{for all}\;|x|>R, \;\hbox{and large}\;n.
$$
 Then it is standard to prove that $\bar{u}_n(x)\rightarrow 0$ exponentially and uniformly in large $n$. Thus, $\displaystyle\lim_{x\rightarrow \infty}\bar{u}(x)=0$.
A similar argument as that in {\sc Step }1 implies that $\bar{u}$ is a weak solution of
$$
   -\De \bar{u}+\bar{u}=\mu_1 \bar{u}^3,\quad \bar{u}(x)\rightarrow 0\;\hbox{as}\;|x|\rightarrow \infty.
$$
So we obtain that $\bar{u}=U_{1,\mu_1}$ and thus $\bar{u}_n(x)\to U_{1,\mu_1}(x)$ in $C_{loc}^{2}(\R^N)$.

b) Using the transformations $\bar{\la}_n:=\frac{1}{\la_n}\to 0$, $\bar{u}_{n}(x):=\frac{1}{\sqrt{\la_n}}v_{n}\left(x/\sqrt{\la_n}\right)$ and $\bar{v}_{n}(x):=\frac{1}{\sqrt{\la_n}}u_{n}\left(x/\sqrt{\la_n}\right)$, we see that $(u_{n},v_{n})$ is a solution to
$$\begin{cases}
-\De u+\la_n u=\mu_1 u^3+\be uv^2\;\quad&\hbox{in}\;\R^N\\
-\De v+v=\mu_2v^3+\be vu^2&\hbox{in}\;\R^N
\end{cases}$$
if and only if $(\bar{u}_{n},\bar{v}_{n})$ is a
solution to
\beq\lab{20190622-we1}
\begin{cases}
-\De u+\bar{\la}_n u=\mu_2 u^3+\be uv^2\;\quad&\hbox{in}\;\R^N,\\
-\De v+v=\mu_1v^3+\be vu^2&\hbox{in}\;\R^N.
\end{cases}
\eeq
We can apply the conclusion of a) to system \eqref{20190622-we1} and obtain that
\[
  \left(\frac{1}{\sqrt{\bar{\la}_n}}\bar{u}_n\left(x/\sqrt{\bar{\la}_n}\right),\bar{v}_n(x)\right) \to (U_{1,\mu_2}(x),U_{1,\mu_1}(x))\quad\hbox{in}\;C_{loc}^2(\R^N)\times C_{loc}^2(\R^N),
\]
that is,
$$
   \left(\frac{1}{\sqrt{\la_n}}u_n\left(x/\sqrt{\la_n}\right),v_n(x)\right) \to (U_{1,\mu_1}(x),U_{1,\mu_2}(x))
     \quad\hbox{in }\;C_{loc}^2(\R^N)\times C_{loc}^2(\R^N).
$$
\ep

\bc\lab{cor:rho}
a) If $(u_{n},v_{n})$ is a positive radial solution to equation \eqref{eq:bif-problem} with $\la=\la_n$ and $\la_n\to 0$ then
$\rho(\la_n,\beta,u_n,v_n)\to +\infty$.

b) If $(u_{n},v_{n})$ is a positive radial solution to equation \eqref{eq:bif-problem} with $\la=\la_n$ and $\la_n\to \infty$ then
$\rho(\la_n,\beta,u_n,v_n)\to 0$.
\ec

\bp
a) Lemma \ref{lem:la-to-0-or-infty} $\bar{u}_n(x):=\frac{1}{\sqrt{\la_n}}u_n(\frac{x}{\sqrt{\la_n}}) \to U_{1,\mu_1}(x)$. So we have that
$$|u_n|_2^2=\la_{n}^{-\frac{1}{2}}|\bar{u}_n|_2^2\to +\infty$$
and
$$|v_n|_2^2\to |U_{1,\mu_2}|_2^2.$$
Hence, $\rho(\la_n,\beta,u_n,v_n)\to +\infty$.

b) Apply a similar argument as in a), and note that $\la_n\to \infty$, we have that
$$|u_n|_2^2=\la_{n}^{-\frac{1}{2}}|\bar{u}_n|_2^2\to 0.$$
\ep

\begin{altproof}{Theorem~\ref{thm:sirakov}}
a) Suppose there exists two families of positive solutions $(u_{\lambda}^{(1)},v_{\lambda}^{(1)})$ and $(u_{\lambda}^{(2)},v_{\lambda}^{(2)})$ to problem \eqref{eq:bif-problem} with $\lambda\rightarrow 0^+$. Let
$$\big(\bar{u}_{\lambda}^{(i)}(x), \bar{v}_{\lambda}^{(i)}(x)\big):=\left(\frac{1}{\sqrt{\lambda}}u_{\lambda}^{(i)}(\frac{x}{\sqrt{\lambda}}), v_{\lambda}^{(i)}(x)\right),\quad i=1,2.$$
Then $\big(\bar{u}_{\lambda}^{(1)}(x), \bar{v}_{\lambda}^{(1)}(x)\big),\big(\bar{u}_{\lambda}^{(2)}(x), \bar{v}_{\lambda}^{(2)}(x)\big)\in E$ are two families of positive solutions to problem
\beq\tag{$P_\lambda$}
\begin{cases}
-\Delta u(x)+u(x)=\mu_1 u(x)^3+\beta u(x)\left(\frac{1}{\sqrt{\lambda}}v(\frac{x}{\sqrt{\lambda}})\right)^2\quad &\hbox{in}\;\R^N,\\
-\Delta v(x)+v(x)=\mu_2v(x)^3+\beta v(x)\left(\sqrt{\lambda}u(\sqrt{\lambda}x)\right)^2&\hbox{in}\;\R^N,\\
0<u,v\in H^1(\R^N), N= 3.
\end{cases}
\eeq
 By Lemma \ref{lem:la-to-0-or-infty},
$$\big(\bar{u}_{\lambda}^{(i)}(x), \bar{v}_{\lambda}^{(i)}(x)\big)\rightarrow (U_{1,\mu_1}, U_{1,\mu_2})\;\hbox{in}\; C_{loc}^{2}(\R^N)\times C_{loc}^{2}(\R^N), i=1,2.$$
Indeed, one can prove that this convergence also holds in $E$ due to the fact that $\bar{u}_{\lambda}^{i}(x)\rightarrow 0$ exponentially and uniformly in small $\lambda$.

\noindent
{\bf Case 1: $\displaystyle\limsup_{\lambda\rightarrow 0^+}\frac{|\bar{v}_{\lambda}^{(1)}-\bar{v}_{\lambda}^{(2)}|_{L^\infty(\R^N)}}{\lambda |\bar{u}_{\lambda}^{(1)}-\bar{u}_{\lambda}^{(2)}|_{L^\infty(\R^N)}}<\infty$}

We study the normalization
$$\xi_\lambda:=\frac{\bar{u}_{\lambda}^{(1)}-\bar{u}_{\lambda}^{(2)}}{|\bar{u}_{\lambda}^{(1)}-\bar{u}_{\lambda}^{(2)}|_{L^\infty(\R^N)}},$$
Then up to a subsequence $\xi_\lambda\rightarrow \xi$ in $C_{loc}^{2}(\R^N)$.
Then we have
\begin{align*}
&\frac{1}{|\bar{u}_{\lambda}^{(1)}-\bar{u}_{\lambda}^{(2)}|_{L^\infty(\R^N)}}
\left[\mu_1\left(\bar{u}_{\lambda}^{(1)}\right)^3-\mu_1\left(\bar{u}_{\lambda}^{(2)}\right)^3\right]\\
=&\mu_1 \xi_\lambda \left[\left(\bar{u}_{\lambda}^{(1)}\right)^2
+\bar{u}_{\lambda}^{(1)}\bar{u}_{\lambda}^{(2)}
+\left(\bar{u}_{\lambda}^{(2)}\right)^2\right]\\
\rightarrow& 3\mu_1 U_{1,\mu_1}^{2} \xi \;\hbox{in}\;C_{loc}^{2}(\R^N)\;\hbox{as}\;\lambda\rightarrow 0,
\end{align*}
and
\begin{align*}
&\frac{1}{|\bar{u}_{\lambda}^{(1)}-\bar{u}_{\lambda}^{(2)}|_{L^\infty(\R^N)}}\left[\beta \bar{u}_{\lambda}^{(1)}(x)\left(\frac{1}{\sqrt{\lambda}}\bar{v}_{\lambda}^{(1)}(\frac{x}{\sqrt{\lambda}})\right)^2-\beta \bar{u}_{\lambda}^{(2)}(x)\left(\frac{1}{\sqrt{\lambda}}\bar{v}_{\lambda}^{(2)}(\frac{x}{\sqrt{\lambda}})\right)^2\right]\\
=&\frac{1}{|\bar{u}_{\lambda}^{(1)}-\bar{u}_{\lambda}^{(2)}|_{L^\infty(\R^N)}}\left[\beta \bar{u}_{\lambda}^{(1)}(x)\left(\frac{1}{\sqrt{\lambda}}\bar{v}_{\lambda}^{(1)}(\frac{x}{\sqrt{\lambda}})\right)^2-\beta \bar{u}_{\lambda}^{(2)}(x)\left(\frac{1}{\sqrt{\lambda}}\bar{v}_{\lambda}^{(1)}(\frac{x}{\sqrt{\lambda}})\right)^2\right]\\
&+\frac{1}{|\bar{u}_{\lambda}^{(1)}-\bar{u}_{\lambda}^{(2)}|_{L^\infty(\R^N)}}\left[\beta \bar{u}_{\lambda}^{(2)}(x)\left(\frac{1}{\sqrt{\lambda}}\bar{v}_{\lambda}^{(1)}(\frac{x}{\sqrt{\lambda}})\right)^2-\beta \bar{u}_{\lambda}^{(2)}(x)\left(\frac{1}{\sqrt{\lambda}}\bar{v}_{\lambda}^{(2)}(\frac{x}{\sqrt{\lambda}})\right)^2\right]\\
=&\beta\xi_\lambda \left(\frac{1}{\sqrt{\lambda}}\bar{v}_{\lambda}^{(1)}(\frac{x}{\sqrt{\lambda}})\right)^2+\beta \bar{u}_{\lambda}^{(2)}(x)\left(\bar{v}_{\lambda}^{(1)}(\frac{x}{\sqrt{\lambda}})+\bar{v}_{\lambda}^{(2)}(\frac{x}{\sqrt{\lambda}})\right)
\frac{\bar{v}_{\lambda}^{(1)}(\frac{x}{\sqrt{\lambda}})-\bar{v}_{\lambda}^{(2)}(\frac{x}{\sqrt{\lambda}})}{\lambda|\bar{u}_{\lambda}^{(1)}-\bar{u}_{\lambda}^{(2)}|_{L^\infty(\R^N)}}.
\end{align*}

For any $h\in H^1(\R^3)$,
one can prove that
\beq\lab{20200112-e1}
\lim_{\lambda\rightarrow 0} \int_{\R^3} \beta\xi_\lambda \left(\frac{1}{\sqrt{\lambda}}\bar{v}_{\lambda}^{(1)}(\frac{x}{\sqrt{\lambda}})\right)^2 h dx=0
\eeq
and
$$
\lim_{\lambda\rightarrow 0} \int_{\R^3}
\bar{u}_{\lambda}^{(2)}(x)\bar{v}_{\lambda}^{(i)}(\frac{x}{\sqrt{\lambda}})h(x)dx=0, i=1,2.
$$
So we see that $\xi$ is a weak solution to
\beq\lab{20200115-e1}
-\Delta \xi+\xi=3\mu_1 U_{1,\mu_1}^{2}\xi.
\eeq
By $|\xi|_{L^\infty}=1$, a standard elliptic estimation indicate that $\xi$ is a strong solution. Then by the decay of $U_{1,\mu_1}$,  applying the comparison principle, we can obtain that $\xi$ exponentially decay to $0$ as $|x|\rightarrow \infty$. Hence, $\xi\in H^1(\R^3)$ and then \eqref{20200115-e1} implies that
$$\xi=\sum_{i=1}^{3} b_i \frac{\partial U_{1,\mu_1}}{\partial x_i}$$
for some suitable $b_i\in \R$.
On the other hand, by the definition, we see that $\xi$ is of radial, and thus $b_i=0, i=1,2,3$. So $\xi=0$, a contradiction. Hence,
$$\bar{u}_{\lambda}^{(1)}\equiv \bar{u}_{\lambda}^{(2)}\;\hbox{for small $\lambda$},$$
and then we also have
$$\bar{v}_{\lambda}^{(1)}\equiv \bar{v}_{\lambda}^{(2)}\;\hbox{for small $\lambda$}$$
due to that
$$\frac{1}{\sqrt{\lambda}}v_{\lambda}^{(i)}(\frac{x}{\sqrt{\lambda}})=\left(\frac{-\Delta \bar{u}_{\lambda}^{(i)}+ \bar{u}_{\lambda}^{(i)}-\mu_1 \left(\bar{u}_{\lambda}^{(i)}\right)^3}{\beta \bar{u}_{\lambda}^{(i)}}\right)^{\frac{1}{2}}, \quad i=1,2.$$

\noindent
{\bf Case 2: $\displaystyle\limsup_{\lambda\rightarrow 0^+}\frac{|\bar{v}_{\lambda}^{(1)}-\bar{v}_{\lambda}^{(2)}|_{L^\infty(\R^N)}}{\lambda |\bar{u}_{\lambda}^{(1)}-\bar{u}_{\lambda}^{(2)}|_{L^\infty(\R^N)}}=\infty$}

In this case, we study the normalization
$$\eta_\lambda:=\frac{\bar{v}_{\lambda}^{(1)}-\bar{v}_{\lambda}^{(2)}}{|\bar{v}_{\lambda}^{(1)}-\bar{v}_{\lambda}^{(2)}|_{L^\infty(\R^N)}},$$
And up to a subsequence,  $\eta_\lambda\rightarrow \eta$ in $C_{loc}^{2}(\R^N)$.
Apply a similar argument as above, we obtain that
$$-\Delta \eta+\eta=3 U_{1,\mu_2}^{2}\eta.$$
By $\eta$ is a radial function,  we also obtain that
$$\bar{v}_{\lambda}^{(1)}\equiv \bar{v}_{\lambda}^{(2)}\;\hbox{for small}\;\lambda,$$
and
$$\bar{u}_{\lambda}^{(1)}\equiv \bar{u}_{\lambda}^{(2)}\;\hbox{for small}\;\lambda$$
by
$$\sqrt{\lambda}\bar{u}_{\lambda}^{(i)}(\sqrt{\lambda}x)=\left(\frac{-\Delta \bar{v}_{\lambda}^{(i)}+\bar{v}_{\lambda}^{(i)}-\mu_2 \left(\bar{v}_{\lambda}^{(i)}\right)^3}{\beta \bar{v}_{\lambda}^{(i)}}\right)^{\frac{1}{2}}, \quad i=1,2.$$

Combining the cases 1 and 2, we see that \eqref{eq:bif-problem} has at most one positive solution for $\lambda$ small enough. And using the transformation in Lemma \ref{lem:la-to-infty}, one can prove the case of $\lambda$ large.

b) It is well known that \eqref{eq:system} has a mountain pass type solution for $\be \le \mu_2\tau_0 < \be_{2}(\la)=\min\{\be_{1}(\la),\be_{2}(\la)\}$ for $\la>0$ small. It follows from a) that this is unique. The second statement in Theorem~\ref{thm:sirakov}~b) for $\be\le\mu_1\tau_0$ follows by applying a transformation as in the proof of Lemma~\ref{lem:la-to-infty}.
\end{altproof}

%%%%%%%%%%%%%%%%%%%%%
\s{Proof of Theorem~\ref{Main-th} and Proposition~\ref{prop:nonexist-1}}\lab{sec:proofs}
%%%%%%%%%%%%%%%%%%%%%
Due to Lemma~\ref{corres-1} it is sufficient to consider the case $\la_1=\la$ and $\la_2=1$, i.e.\ system~\eqref{eq:bif-problem}. \\

\begin{altproof}{Theorem~\ref{Main-th}}
a) For $\be \le \tau_0\min\{\mu_1,\mu_2\}$ the existence of normalized solutions for every $a,b>0$ follows from Corollaries~\ref{cor:beta-small} and \ref{cor:rho}. For $\be \ge \tau_0\max\{\mu_1,\mu_2\}$ let $\cs_i^\be$, $i=1,2$, be the connected sets of positive solutions from Proposition~\ref{prop:global}~e). If $\cs_1^\be \cap \cs_2^\be \ne \emptyset$ then the existence of normalized solutions for every $a,b>0$ follows from Corollary~\ref{cor:cs1=cs2}. Now we suppose $\cs_1^\be \cap\cs_2^\be = \emptyset$. Then Proposition~\ref{prop:global}~e) yields that $P_1(\cs_i^\be)$ contains one of the intervals $(0,\ell_i(\be))$ or $(\ell_i(\be),\infty)$, $i=1,2$. If $(\ell_1(\be),\infty) \subset P_1(\cs_1^\be)$ then the existence of normalized solutions for every $a,b>0$ follows from Corollary~\ref{cor:rho}. The same argument applies if $(0,\ell_2(\be)) \subset P_1(\cs_2^\be)$. Now we show that the case $\cs_1^\be \cap\cs_2^\be = \emptyset$ and $(0,\ell_2(\be)) \not\subset P_1(\cs_2^\be)$ cannot happen, concluding the proof of a). Similarly one can show that $\cs_1^\be \cap\cs_2^\be = \emptyset$ and $(\ell_1(\be),\infty) \not\subset P_1(\cs_1^\be)$ leads to a contradiction.

Suppose by contradiction that $\cs_1^\be \cap \cs_2^\be = \emptyset$ and $(0,\ell_2(\be)) \not\subset P_1(\cs_2^\be)$. Then $(\ell_2(\be),\infty)\subset P_1(\cs_2^\be)$. Recall from Theorem~\ref{thm:sirakov}~a) that \eqref{eq:bif-problem} has at most one solution for $\la$ large. It follows that there exists a family $(\la,\be,u_{\la,\be},v_{\la,\be})\in X$, $\la\ge\tilde{\la}(\be)$, so that
\[
  \cs^\be\cap\big([\tilde{\la}(\be),\infty)\times\P\big) = \cs_1^\be\cap\big([\tilde{\la}(\be),\infty)\times\P\big)
       = \{(\la,\be,u_{\la,\be},v_{\la,\be}):\la\ge\tilde{\la}(\be)\}.
\]
The fixed point index computations in Section~\ref{sec:global}, in particular Propositions~\ref{prop:bif}, \ref{prop:i_infty} and Corollary~\ref{cor:beta_i}, imply for $\la\ge\tilde{\la}(\be)$:
\beq\lab{eq:ind-la-large}
\begin{aligned}
  \ind_\P\big(\A_{\la,\be},(u_{\la,\be},v_{\la,\be})\big)
   &= i_\infty - \ind_{\P}\big(\A_{\la,\be},(U_{\la,\mu_1},0)\big)\\
   &\hspace{.5cm} - \ind_{\P}\big(\A_{\la,\be},(0,U_{1,\mu_2})\big)
        - \ind_{\P}\big(\A_{\la,\be},(0,0)\big)\\
    &= 0+0+1-1 = 0
\end{aligned}
\eeq
Observe that $\ct_2^\be\cup\cs_2^\be$ is a connected component of the set
$\cz=\ct_0\cup\ct_1\cup\ct_2\cup\cs$
of all solutions because $\cs_1^\be \cap \cs_2^\be = \emptyset$. Then there exists an open set $\co \subset X^\be$ with the following properties:
\begin{itemize}
\item[(i)]$\ct_2^\be\cup\cs_2^\be \subset \co$
\item[(ii)] $\cz\cap\pa \co=\emptyset$
\item[(iii)] There exists $\de>0$ so that
    \[
      \co\cap \big((0,\de]\times\{\be\}\times\P\big) = \big\{(\la,\be,u,v):\la\in(0,\de],\ (u,v)\in B_\de(0,U_{1,\mu_2})\big\}
    \]
\end{itemize}
The last property (iii) can be achieved because  $(0,\ell_2(\be)) \not\subset P_1(\cs_2^\be)$, hence $\cs_2^\be\subset [\de,\infty)\times\{\be\}\times\P$ for some small $\de>0$. Using the notation $\co_{\la,\be}:=\{(u,v)\in\P: (\la,\be,u,v)\in \co\}$ it follows for $\la\ge\tilde{\la}(\be)$ that:
\[
\begin{aligned}
  \ind_\P\big(\A_{\la,\be},(u_{\la,\be},v_{\la,\be})\big)
       &= \ind_\P\big(\A_{\la,\be},\co_{\la,\be}\big) - \ind_{\P}\big(\A_{\la,\be},(0,U_{1,\mu_2})\big)\\
       &= \ind_\P\big(\A_{\de,\be},\co_{\de,\be}\big) - \ind_{\P}\big(\A_{\la,\be},(0,U_{1,\mu_2})\big)\\
       &= \ind_\P\big(\A_{\de,\be},(0,U_{1,\mu_2})\big) - \ind_{\P}\big(\A_{\la,\be},(0,U_{1,\mu_2})\big)\\
       &= 0+1 = 1
\end{aligned}
\]
This contradicts \eqref{eq:ind-la-large}.

b) We only prove the case $\mu_2<\mu_1$. The case $\mu_1<\mu_2$ can then be deduced using the transformation from the proof of Lemma~\ref{lem:la-to-infty}. Let $\cs_2^\be$ be the connected set of positive solutions from Proposition~\ref{prop:global}~d). Then Proposition~\ref{prop:global}~d) yields that $P_1(\cs_2^\be)$ contains one of the intervals $(0,\ell_2(\be))$ or $(\ell_2(\be),\infty)$. If $(0,\ell_2(\be)) \subset P_1(\cs_2^\be)$ then the existence of normalized solutions for every $a,b>0$ follows from Corollary~\ref{cor:rho}. If $(\ell_2(\be),\infty)\subset P_1(\cs_2^\be)$ then
\[
  \de:=\max_{(\la,\be,u,v)\in \cs_2^\be} \rho(\la,\be,u,v) > 0.
\]
Since $\rho(\la,\be,u,v) \to 0$ as $\la\to\infty$, and as $\la\to\ell_2(\be)$ on $\cs_2^\be$, we see that $\rho(\cs)\supset (0,\de]$.

Finally, if $\be \in (\tau_0\mu_2,\mu_2)$ then there exists the solution $(1,\be,u_\be,v_\be)\in\cs$ from Remark~\ref{rem:unique}, which has fixed point index 1. Let $\cs_0^\be\subset\cs^\be$ be the connected component of $(1,\be,u_\be,v_\be)$ in $\cs^\be$. An index count as above yields that $P_1(\cs_0^\be)\subset\R^+$ is bounded away from $0$. Since it cannot bifurcate from $\ct_1$ it must bifurcate from $\ct_2$, i.e.\ $\cs_3^\be=\cs_2^\be$. This implies
\[
  \de\ge \rho(1,\be,u_\be,v_\be) = \sqrt{\frac{\be-\min\{\mu_1,\mu_2\}}{\be-\max\{\mu_1,\mu_2\}}}.
\]
\end{altproof}

\begin{altproof}{Proposition~\ref{prop:nonexist-1}}
We only prove the case of $\mu_2\leq \beta\leq \tau \mu_1$, the second part result is easy by using the transformation from the proof of Lemma~\ref{lem:la-to-infty}.  By Theorem \ref{thm:nonexist-2} b), there exists $\eta_2(\beta)>0$ such that problem \eqref{eq:bif-problem} has no positive solution provided $\la<\eta_2(\beta)$.  On the other hand, by Theorem \ref{thm:sirakov} b), problem \eqref{eq:bif-problem} has a unique positive solution $(u_\la,v_\la)$, which is of mountain pass type, for $\la \ge \tilde{\la}(\be)$ large enough. By Corollary \ref{cor:rho}, we have that $\rho(\la,\beta,u_\la,v_\la)\to 0$ as $\la\to \infty$.
So
$$q_1:=\{\rho(\la,\beta,u_\la,v_\la),\la \geq \tilde{\la}(\beta)\;\}<\infty.$$
Observe that according to Proposition \ref{prop:bif} a), see also \cite[Lemma 2.1]{BartschWangWei2007},
$$\sup_{(\la,\be,u,v)\in \cs^\be, \eta_2(\beta)\leq\la \leq\tilde{\la}(\be)} \left(|u|_2^2+|v|_2^2\right)<\infty.$$
Then we have that
$$
   q_2:=\sup \{\rho(\la,\beta,u,v):(\la,\beta,u,v)\in \cs^\be,\ \eta_2(\beta)\leq\la \leq\tilde{\la}_\beta\}<\infty.
$$
Indeed, if there exists a sequence $(\la_n,\be,u_n,v_n)$ with $\la_n\to \la\in [\eta_2(\beta),\tilde{\la}_\beta]$ such that $\rho(\la_n,\be,u_n,v_n)\rightarrow \infty$. Then we see that $|v_{n}|_2^2\to 0$ and it is standard to prove that $(u_{n}, v_{n})\to (U_{\la,\mu_1}, 0)$ in $H^1(\R^N)$. And thus, $\displaystyle\beta=\beta_{1}(\la)>\lim_{\la\to \infty} \beta_{1}(\la)=\tau_0 \mu_1$, a contradiction. Then $q:=\max\{q_1,q_2\}$ is the required bound.

\end{altproof}

\noindent
{\bf Acknowledgements}\\
The author Xuexiu Zhong thanks Zhijie Chen for the valuable discussions when preparing the paper.

%%%%%%%%%%%%%%%%%%%%%%%%%%%%%%%%%%%%%%%%%%%%%%%%%%%%%%%%%%%%%%%%%%%%%%%%%%%%%%

%\bibliographystyle{abbrv}
%\bibliography{ref_normalized-solution}

\vskip0.26in
\begin{thebibliography}{10}

\bibitem{AkhmedievAnkiewicz.1999}
N.~Akhmediev and A.~Ankiewicz:
\newblock Partially coherent solitons on a finite background.
 {\it Phys. Rev. Lett.}, {\bf 82}(13), 2661, 1999.

\bibitem{AlexanderAntman1981}
J.~C.~ Alexander and S.~Antman:
\newblock Global and local behavior of bifurcating multidimensional continua of solutions for multiparameter nonlinear eigenvalue problems.
{\it Arch. Ration. Mech. Anal.}, {\bf 76}(4), 339--354, 1981.

\bibitem{Ambrosetti2006}
A.~Ambrosetti and E.~Colorado:
\newblock Bound and ground states of coupled nonlinear {S}chr\"odinger equations.
 {\it C. R. Math. Acad. Sci. Paris}, {\bf 342}(7), 453--458, 2006.

\bibitem{Ambrosetti2007}
A.~Ambrosetti and E.~Colorado:
\newblock Standing waves of some coupled nonlinear {S}chr\"odinger equations.
 {\it J. Lond. Math. Soc. (2)}, {\bf 75}(1), 67--82, 2007.

\bibitem{BartschDancer2009}
T.~Bartsch and N.~Dancer:
\newblock Poincar\'e-Hopf type formulas on convex sets of Banach spaces.
{\it Topol. Methods Nonlinear Anal.}, {\bf 34}, 213--229, 2009.

\bibitem{Bartsch2010}
T.~Bartsch, N.~Dancer and Z.-Q.~Wang:
\newblock A Liouville theorem, a priori bounds, and bifurcating branches of positive solutions for a nonlinear elliptic system.
{\it Calc. Var. Partial Diifer. Equ.}, {\bf 37}, 345-361, 2010.

\bibitem{Bartsch2018}
T.~Bartsch and L.~Jeanjean:
\newblock Normalized solutions for nonlinear {S}chr\"odinger systems.
 {\it Proc. Roy. Soc. Edinburgh Sect. A}, {\bf 148}(2), 225--242, 2018.

\bibitem{Bartsch2016}
T.~Bartsch, L.~Jeanjean, and N.~Soave:
\newblock Normalized solutions for a system of coupled cubic {S}chr\"odinger
  equations on {$\Bbb{R}^3$}.
{\it J. Math. Pures Appl. (9)}, {\bf 106}(4), 583--614, 2016.

\bibitem{Bartsch2017}
T.~Bartsch and N.~Soave:
\newblock A natural constraint approach to normalized solutions of nonlinear
  {S}chr\"odinger equations and systems.
{\it J. Funct. Anal.}, {\bf 272}(12), 4998--5037, 2017.

\bibitem{Bartsch2019}
T.~Bartsch and N.~Soave:
\newblock Multiple normalized solutions for a competing system of Schr\"odinger equations.
 {\it Calc. Var.}, {\bf 58:22}, 2019.

\bibitem{BartschWang2006}
T.~Bartsch and Z.-Q.~Wang:
\newblock Note on ground states of nonlinear Schr\"odinger systems.
{\it J. Partial Differ. Equ.}, {\bf 19}, 200-207, 2006.

\bibitem{BartschWangWei2007}
T.~Bartsch, Z.-Q.~Wang and J.~Wei:
\newblock Bound states for a coupled Schr\"odinger system:
{\it J. Fixed Point Theory Appl.}, {\bf 2}(2), 353--367, 2007.

\bibitem{BuscaSirakov2001}
J.~Busca and B.~Sirakov:
\newblock Symmetry results for semilinear elliptic systems in the whole space.
{\it J. Differential Equations}, {\bf 163}(1), 41--56, 2001.

\bibitem{Chen2013a}
Z.~Chen and W.~Zou:
\newblock An optimal constant for the existence of least energy solutions of a
  coupled {S}chr\"odinger system.
 {\it Calc. Var. Partial Differential Equations}, {\bf 48}(3-4), 695--711, 2013.

\bibitem{CrandallRabinowitz1971}
M.~Crandall and P.~Rabinowitz:
\newblock Bifurcation form fimple eigenvalues.
{\it J. Funct. Anal.}, {\bf 8}, 321-340, 1971.

\bibitem{CrandallRabinowitz1973}
M.~Crandall and P.~Rabinowitz:
\newblock Bifurcation, perturbation of simple eigenvalues and linearized stability.
{\it Arch. Ration. Mech. Anal.}, {\bf 52}, 161-180, 1973.

\bibitem{DancerWei2009}
E.~Dancer and J.~Wei:
\newblock Spike solutions in coupled nonlinear Schr\"odinger equations with attractive interaction.
{\it Trans. Amer. Math. Soc.},{\bf 361}, 1189-1208, 2009.

\bibitem{Esry1998}
B.~D.~Esry, C.~H.~Greene, J.~P.~Burke Jr., and J.~L.~Bohn:
\newblock Hartree-fock theory for double condensates.
{\it Phys. Rev. Lett.}, {\bf 78} (3594), 1997.

\bibitem{Frantzeskakis2010}
D.~J. Frantzeskakis:
\newblock Dark solitons in atomic Bose-Einstein condensates: from theory to experiments.
 {\it J. Phys. A: Math. Theor.}, {\bf 43}, 2010.

\bibitem{GidasSpruck.1981}
B.~Gidas and J.~Spruck:
\newblock Global and local behavior of positive solutions of nonlinear elliptic equations.
 {\it Comm. Pure Appl. Math.}, {\bf 34}(4), 525--598,
  1981.

\bibitem{GouJeanjean2018}
T.~Gou and L.~Jeanjean:
\newblock Multiple positive normalized solutions for nonlinear Schr\"odinger systems.
{\it Nonlinearity} {\bf 31}(2), 2319--2345, 2018.

\bibitem{Ikoma2011}
N.~Ikoma and K.~Tanaka:
\newblock A local mountain pass type result for a system of nonlinear Schr{\"o}dinger equations.
{\it Calc. Var. Partial Differential Equations},
  {\bf 40}(3-4), 449--480, 2011.

\bibitem{Kwong.1989}
M.~K. Kwong:
\newblock Uniqueness of positive solutions of $\De u- u+ u^p= 0$ in $\R^n$.
 {\it Arch. Ration. Mech. Anal.}, {\bf 105}(3), 243--266, 1989.

\bibitem{LinWei.2005}
T.-C. Lin and J.~Wei:
\newblock Ground state of n coupled nonlinear Schr{\"o}dinger equations in $\R^n, n\leq 3$.
 {\it Comm. Math. Phys.}, {\bf 255}(3) ,629--653, 2005.

\bibitem{MaiaMontefuscoPellacci.2006}
L.~Maia, E.~Montefusco, and B.~Pellacci:
\newblock Positive solutions for a weakly coupled nonlinear Schr{\"o}dinger
  system.
 {\it J. Differential Equations}, {\bf 229}(2), 743--767, 2006.

\bibitem{Mandel2016}
R.~Mandel:
\newblock Minimal energy solutions and infinitely many bifurcating branches for
  a class of saturated nonlinear {S}chr\"odinger systems.
{\it Adv. Nonlinear Stud.}, {\bf 16}(1), 95--113, 2016.

%\bibitem{NiTakagi1993}
%W.~Ni and I.~Takagi:
%\newblock Locating the peaks of least energy solutions to a semilinear Neumann problem.
 %{\it Duke Math. J.}, {\bf 70}(2), 1--20, 1993.

\bibitem{Noris2015}
B.~Noris, H.~Tavares, and G.~Verzini:
\newblock Stable solitary waves with prescribed $L^2$-mass for the cubic Schr\"odinger system with trapping potential.
 {\it Discrete Contin. Dyn. Syst.}, {\bf 35}(12), 6085--6112, 2015.

\bibitem{Noris2019}
B.~Noris, H.~Tavares, and G.~Verzini:
\newblock Normalized solutions for nonlinear Schr\"odinger systems on bounded
  domains.
 {\it Nonlinearity}, {\bf 32}(3), 1044--1072, 2019.

\bibitem{Sirakov2007}
B.~Sirakov:
\newblock Least energy solitary waves for a system of nonlinear Schr{\"o}dinger
  equations in $\R^n$.
 {\it Comm. Math. Phys.}, {\bf 271}(1), 199--221, 2007.

\bibitem{Soave2015}
N.~Soave:
\newblock On existence and phase separation of solitary waves for nonlinear
  Schr\"odinger systems modelling simultaneous cooperation and competition.
 {\it Calc. Var. Partial Differential Equations}, {\bf 53}(3-4), 689--718, 2015.

\bibitem{Soave2016}
N.~Soave and H.~Tavares:
\newblock New existence and symmetry results for least energy positive
  solutions of {S}chr\"odinger systems with mixed competition and cooperation
  terms.
{\it J. Differential Equations}, {\bf 261}(1), 505--537, 2016.

\bibitem{Terracini2009}
S.~Terracini and G.~Verzini:
\newblock Multipulse phases in {$k$}-mixtures of {B}ose-{E}instein condensates.
 {\it Arch. Ration. Mech. Anal.}, {\bf 194}(3), 717--741, 2009.

\bibitem{Timmermans1998}
E.~Timmermans:
\newblock Phase separation of Bose-Einstein condensates.
 {\it Phys. Rev. Lett.}, {\bf 81}, 5718--5721, 1998.

\bibitem{Wei2008}
J.~Wei and T.~Weth:
\newblock Radial solutions and phase separation in a system of two coupled
  {S}chr\"odinger equations.
 {\it Arch. Ration. Mech. Anal.}, {\bf 190}(1), 83--106, 2008.

\bibitem{Wei2012}
J.~Wei and W.~Yao:
\newblock Uniqueness of positive solutions to some coupled nonlinear Schr\"odinger equations.
{\it Commun. Pure. Appl. Anal.}, {\bf 11}, 1003-1011, 2012.

\end{thebibliography}

%%%%%%%%%%%%%%%%%%%%%%%%%%%%%%%%%%%%%%%%%%%%%%%%%%%%%%%%%%%%%%%%%
\end{CJK*}
\end{document}